\documentclass[english]{elsarticle}
\usepackage[T1]{fontenc}
\usepackage[utf8]{inputenc}
\usepackage{babel}
\usepackage{verbatim}
\usepackage{amsmath}
\usepackage{amsthm}
\usepackage{amssymb}
\usepackage{graphicx}
\usepackage[unicode=true,
 bookmarks=true,bookmarksnumbered=false,bookmarksopen=false,
 breaklinks=false,pdfborder={0 0 1},backref=false,colorlinks=false]
 {hyperref}

\makeatletter
\theoremstyle{plain}
\newtheorem{thm}{\protect\theoremname}
  \theoremstyle{remark}
  \newtheorem{rem}[thm]{\protect\remarkname}
  \theoremstyle{definition}
  \newtheorem{example}[thm]{\protect\examplename}

\usepackage{tikz}

\makeatother

  \providecommand{\examplename}{Example}
  \providecommand{\remarkname}{Remark}
\providecommand{\theoremname}{Theorem}

\begin{document}

\begin{frontmatter}{}

\title{High-order implicit palindromic discontinuous Galerkin method for
kinetic-relaxation approximation}

\tnotetext[t2]{This work has been supported by Eurofusion, ANR project ``EXAMAG'' and BPI France project ``HOROCH''.}

\author{David Coulette, Emmanuel Franck, Philippe Helluy, Michel Mehrenberger,
Laurent Navoret}

\ead{helluy@unistra.fr}

\address{IRMA, Univ. Strasbourg, 7 rue Descartes, Strasbourg, France \& Inria
TONUS}

\begin{abstract}
We construct a high order discontinuous Galerkin method for solving
general hyperbolic systems of conservation laws. The method is CFL-less,
matrix-free, has the complexity of an explicit scheme and can be of
arbitrary order in space and time. The construction is based on: (a)
the representation of the system of conservation laws by a kinetic
vectorial representation with a stiff relaxation term; (b) a matrix-free,
CFL-less implicit discontinuous Galerkin transport solver; and (c)
a stiffly accurate composition method for time integration. The method
is validated on several one-dimensional test cases. It is then applied
on two-dimensional and three-dimensional test cases: flow past a cylinder,
magnetohydrodynamics and multifluid sedimentation.
\end{abstract}

\end{frontmatter}{}

\tableofcontents{}\global\long\def\vf{\mathbf{{f}}}
\global\long\def\vV{\mathbf{V}}
\global\long\def\vx{\mathbf{x}}
\global\long\def\v#1{\mathbf{#1}}
\global\long\def\vq{\mathbf{q}}
\global\long\def\vw{\mathbf{w}}
\global\long\def\ddim{D}
\global\long\def\dorder{d}
\global\long\def\normal{\mathbf{n}}
\global\long\def\flux{\mathbf{q}}
\global\long\def\Vu{\mathbf{u}}
\global\long\def\Vvi{\mathbf{v}_{i}}
\global\long\def\vg{\mathbf{g}}
\global\long\def\vs{\mathbf{s}}
\global\long\def\vtf{\mathbf{{\tilde{f}}}}
\global\long\def\Skd{{\sum_{k=1}^{D}}}
\global\long\def\Sjd{{\sum_{j=1}^{D}}}
\global\long\def\vQ{\mathbf{Q}}
\global\long\def\vP{\mathbf{P}}
\global\long\def\Id{\text{\text{Id}}}
\global\long\def\vL{\mathbf{L}}
\global\long\def\vNh{\mathbf{\mathbf{N}}_{h}^{\tau}}
\global\long\def\vN{\mathbf{N}^{\tau}}
\global\long\def\vR{\mathbf{R}^{\tau}}
\global\long\def\vfh{\mathbf{\mathbf{f}}_{h}}
\global\long\def\vM{\mathbf{M}}
\global\long\def\vRz{\mathbf{R}^{0}}
\global\long\def\vSh{\mathbf{S}_{h}}
\global\long\def\vS{\mathbf{S}}
\global\long\def\vT{\mathbf{T}}
\global\long\def\vG{\mathbf{G}}
\global\long\def\vH{\mathbf{H}}

\section{Introduction}

Systems of conservation laws are important mathematical tools for
modelling many phenomena in physics or engineering.

In several practical applications, only some time scales of the model
are interesting and one would like to filter out the smallest time
scales. Classical explicit methods require very small time steps,
because of the CFL stability condition. A standard way to treat the
various time scales is to use an implicit time-stepping scheme. Schemes
of this kind, however, are quite challenging from a computational
point of view: they require inverting large non-linear systems, which
induce high computational and storage costs.

In this paper, we propose an alternative method for solving systems
of conservation laws for a large range of time scales on complex geometries.
The time-marching procedure is CFL-less, but keeps the complexity
of an explicit scheme. In addition, we are able to achieve high order
in space and time.

Our method is based on a vectorial kinetic relaxation scheme described
in \cite{aregba2000discretekinetic,bouchut2004stability,perthame1990boltzmann,natalini1998kinetic}.
The vectorial kinetic scheme is a generalization of the relaxation
scheme of Jin and Xin \cite{jin1995relaxation}. The original system
of conservation laws is replaced by an equivalent kinetic system made
of a small set of transport equations, coupled through a stiff relaxation
source term. The stiffness is measured with a small relaxation time
$\tau>0$. The original system of conservation laws is equivalent
to its kinetic representation in the limit $\tau\to0$.

Many approaches have been proposed in the literature for solving such
kinetic models. It is generally approximated by a splitting method,
in which the transport and relaxation steps are treated separately.
A simple and natural choice for solving the stiff relaxation step
is then to apply a first order implicit scheme in order to avoid instabilities.
On the other hand, the transport step can be solved with several different
methods: explicit upwind schemes \cite{brenier1984averaged,coron1991numerical,aregba2000discretekinetic};
exact characteristic schemes (which are at the base of the Lattice
Boltzmann Method \cite{qian1992lattice,chen1998lattice,he1997lattice});
but also finite volume, finite difference or discontinuous Galerkin
methods (see for instance \cite{nannelli1992lattice,peng1998lattice,shi2003discontinuous,mei1998finite}). 

From a computational point of view, the splitting approach has several
advantages: the transport equations are uncoupled, linear and can
be solved with efficient parallel solvers; the relaxation step is
also embarrassingly parallel and requires only to considering local
ordinary differential equations.

However, the first order splitting introduces too much numerical diffusion
for practical applications. Therefore many works have been devoted
to the construction of higher order schemes based on improved splitting
approaches.

In the Lattice Boltzmann Method, the accuracy is improved if the relaxation
step is solved with a Crank-Nicolson scheme \cite{dellar2013interpretation}.
The transport and relaxation steps are then interlaced with a Strang
procedure. Because of the stiff relaxation, it has been observed that
without special care in the scheme design, one can observe order reduction
when the relaxation time $\tau\to0$ \cite{jin1995runge}. It is,
however, possible to construct high order Runge-Kutta schemes, mixing
implicit and explicit steps \cite{pareschi2005implicit} that preserve
the accuracy when $\tau\to0$ (Asymptotic Preserving property \cite{jin1999efficient}).

In all the above approaches, because the transport step is solved
by an explicit scheme the whole procedure is still constrained by
a CFL condition on the time step.

The first fundamental aspect of our method is to apply an \textit{implicit}
Discontinuous Galerkin (DG) method instead of an explicit one for
solving the transport equations. In this way, we obtain unconditionally
stable schemes and get rid of the CFL condition. The implicit solver
has almost no additional cost compared to the explicit one. Indeed,
with an upwind numerical flux, the linear system of the implicit DG
method is triangular and, in the end, can be solved explicitly. This
kind of ideas is mentioned in several works. See for instance in \cite{bey1997downwind,wang1999crosswind,coquel2008large,natvig2008fast,moustafa20153d}.
In a recent work we have evaluated the parallel scalability of the
triangular solver \cite{badwaik2017task}.

The second fundamental aspect of our method is the construction of
a \textit{symmetric-in-time} integrator that remains second order
accurate even for vanishing relaxation time $\tau$ (AP property).
The construction is based on a modified Crank-Nicolson procedure and
on essential reversibility properties of the transport equation. Once
a symmetric-in-time integrator is available, it is then very easy
to construct arbitrary order methods with the composition method \cite{suzuki1990fractal,kahan1997composition,mclachlan2002splitting,hairer2006geometric}.
We apply this method for achieving fourth and sixth order time integration
even for vanishing relaxation time $\tau$.

The objective of this paper is first to present the whole construction
of the Palindromic Discontinuous Galerkin Method. Then we will establish
some rigorous properties of the scheme in the simplified linear case.
We will validate the approach on several one-dimensional test cases.
Finally, we will apply it in 2D and 3D for computing Von Karmann streets
and multi-fluid instabilities.

\section{Kinetic relaxation approximation}

We consider a system of $m$ conservation laws in $D$ space dimension.
The unknown $\v w(\v x,t)\in\mathbb{R}^{m}$, depending on space $\vx=(x^{1},\ldots,x^{\ddim})\in\mathbb{R}^{\ddim}$
and time $t$> 0, satisfies the following system
\begin{equation}
\partial_{t}\vw+\sum_{k=1}^{\ddim}\partial_{k}\left(\flux^{k}(\vw)\right)=\vs(\vw),\label{eq:macro_limit_system-1}
\end{equation}

where $\flux^{k}(\vw)\in\mathbb{R}^{m}$ are the fluxes in the $k$-th
spatial direction with $1\leq k\leq D$ and $s(\vw)\in\mathbb{R}^{m}$
is a general source term. For any function $g(\vx,t)$, $\partial_{k}g(\vx,t)$
stands for the partial derivative of $g$ with respect to $x^{k}$.

The kinetic BGK representation aims at considering (\ref{eq:macro_limit_system-1})
as a singular limit of a linear kinetic equation with a source term.

\subsection{Kinetic BGK equation}

The macroscopic quantity $\v w(\v x,t)\in\mathbb{R}^{m}$ is associated
to a vectorial distribution function $\vf(\vx,t)\in\mathbb{R}^{n_{v}}$,
with $n_{v}>m$, through a linear transformation
\begin{equation}
\v w=\v P\vf,\label{eq:micro_to_macro}
\end{equation}
where $\v P$ is a constant $m\times n_{v}$ matrix. Each component
of $\vf(\vx,t)$ corresponds to a discrete velocity, denoted $\v v_{i}=(v_{i}^{1},\ldots,v_{i}^{D})\in\mathbb{R}^{D}$
for $1\leq i\leq n_{v}$. This distribution function satisfies the
following (kinetic) equation
\begin{equation}
\partial_{t}\vf+\Skd\vV^{k}\partial_{k}\vf=\frac{1}{\tau}\Big(\vf^{\text{eq}}\big(\v P\vf\big)-\vf\Big)+\vg(\vf).\label{eq:kin_bgk}
\end{equation}
where, for all $1\leq k\leq\ddim$, $\vV^{k}\in M_{n_{v}}(\mathbb{R})$
are diagonal matrices composed of the $k$-th components of the discrete
velocities
\[
\vV^{k}=\left(\begin{array}{cccc}
v_{1}^{k}\\
 & v_{2}^{k}\\
 &  & \ddots\\
 &  &  & v_{n_{v}}^{k}
\end{array}\right).
\]
The right-hand side of (\ref{eq:kin_bgk}) involves a generic source
term $\vg(\vf)$ and a BGK relaxation term: $\v{N^{\tau}\vf}=(\vf^{\text{eq}}(\v P\vf)-\vf)/\tau$,
where $\tau\ll1$ is a small parameter. In other words, equation (\ref{eq:kin_bgk})
is a coupling of $n_{v}$ transport equations at constant velocities. 

The relaxation term is devised so that the macroscopic quantity $\v w=\v P\vf$
converges to the solution to equation (\ref{eq:macro_limit_system-1})
as $\tau\to0$. To ensure such behaviour, the equilibrium distribution
$\vf^{\text{eq}}(\v P\vf)$ only depends on the macroscopic quantity and
satisfy:
\begin{equation}
\v w=\v P\vf^{\text{eq}}(\v w).\label{eq:Pfeq_w}
\end{equation}

Consequently, multiplying equation (\ref{eq:kin_bgk}) by $\v P$
makes the singular relaxation term vanish and we get
\begin{equation}
\partial_{t}\v P\vf+\Skd\partial_{k}\big(\v P\vV^{k}\vf\big)=\v P\vg(\vf).\label{eq:kin_bgk-1}
\end{equation}

Since $\vf$ formally tends to $\vf^{\text{eq}}(\v w)$ as $\tau\to0$, we
recover the system of conservation (\ref{eq:macro_limit_system-1})
in the limit provided that we have the following relation
\begin{align}
\flux^{k}(\vw) & =\v P\vV^{k}\vf^{\text{eq}}(\v w),\label{eq:flux_kin}\\
\vs(\vw) & =\mathbf{P}\vg(\vf^{\text{eq}}(\v w)).\label{eq:source_kin}
\end{align}
In the next section, we will give examples of such constructions.

As shown in \cite{aregba2000discretekinetic}, at the first order
in $\tau$, the kinetic relaxation system is consistent with
\begin{multline}
\partial_{t}\vw+\Skd\partial_{k}(\vq^{k}(\vw))=\vs+\tau\Skd\Sjd\partial_{k}[\mathcal{D}^{kj}(\vw)\partial_{j}\vw]\\
+\tau\left[\Skd\partial_{k}\mathbf{\Big(P}\vV^{k}\big[\nabla_{\v w}\vf^{\text{eq}}(\vw)\vs(\vw)-\vg(\vf^{\text{eq}}(\v w))\big]\Big)\right]+\mathcal{O}(\tau^{2}),\label{eq:macro_asymptic_first_order}
\end{multline}

where the diffusion tensor $\mathcal{D}$ is defined by

\begin{equation}
\mathcal{D}^{kj}=\mathbf{P}\vV^{k}\vV^{j}\nabla_{\v w}\vf^{\text{eq}}-\nabla_{\v w}\vq^{k}\nabla_{\v w}\vq^{j}.\label{eq:diffusion_vecto_general}
\end{equation}

For the sake of completeness, the proof of this estimate is provided
in Appendix \ref{subsec:Second-order-approximation}.

In view of (\ref{eq:macro_asymptic_first_order}), it is particularly
interesting to take the kinetic source equal to 
\begin{equation}
\vg(\vf)=\nabla_{\v w}\vf^{\text{eq}}(\vP\vf)\,\vs(\vP\vf).\label{eq:source_jacfeq_s}
\end{equation}
This choice directly guaranties the consistency of the source (\ref{eq:source_kin})
since (\ref{eq:Pfeq_w}) implies that $\mathbf{P}\nabla_{\vw}\vf^{\text{eq}}$
reduces to the identity matrix on $\mathbb{R}^{m}$. This particular
form of the source term is actually equivalent to only make
the macroscopic part of the distribution function  evolve (see Remark \ref{rem:source}
below).

The stability of kinetic relaxation models is discussed in \cite{aregba2000discretekinetic,bouchut2004stability,chen1994hyperbolic}.
The mere dissipation of the $L^{2}$ norm, which requires the symmetric
part of the diffusion tensor to be positive, is not sufficient for
nonlinear hyperbolic systems. For such systems, a more appropriate
criterion is the dissipation of an entropy: in its strongest form,
it requires the existence of a strictly convex entropy for the kinetic
system. A weaker requirement is the dissipation of a macroscopic entropy
by the approximated system at the first order in the Chapman-Enskog
expansion. We assume the existence of a convex entropy-flux pair $(\eta(\vw),\vQ^{k}(\vw))$
for (\ref{eq:macro_limit_system-1}) and note $\nabla_{\mathbf{w}}^{2}\eta$
the Hessian matrix of the entropy. From \eqref{eq:macro_asymptic_first_order},
we have
\begin{multline}
\partial_{t}\eta(\mathbf{w)}+\Skd\partial_{k}\big(\vQ^{k}(\mathbf{w)\big)}+\nabla_{\vw}\eta^{T}\vs-\tau\Skd\Sjd\partial_{k}[(\nabla_{\mathbf{w}}\eta(\mathbf{w)})^{T}\mathcal{D}^{kj}(\mathbf{w)}\partial_{j}\vw]\\
=-\tau\Skd\Sjd(\partial_{k}\vw)^{T}(\nabla_{\mathbf{w}}^{2}\eta(\mathbf{w)})^{T}D^{kj}\mathbf{(w)}\partial_{j}\vw,\label{eq:entropy}
\end{multline}
whose r.h.s is dissipative provided the tensor $(\nabla_{\mathbf{w}}^{2}\eta)^{T}\mathcal{D}^{kj}$
is definite non-negative.
\begin{rem}
\label{rem:source}The kinetic source term (\ref{eq:source_jacfeq_s})
makes  the macroscopic variable $\vw=\vP\vf$ evolve according to the
macroscopic source dynamics but leaves the out-of-equilibrium part
$\vtf=\vf-\vf^{\text{eq}}(\vP\vf)$ unchanged. Indeed, considering the differential
equation
\[
\partial_{t}\vf=\nabla_{\mathbf{w}}\vf^{\text{eq}}(\vP\vf)\,\vs(\vP\vf),
\]
we easily show that $\vP\vf$ satisfies the differential equation
\[
\partial_{t}(\vP\vf)=\vs(\vP\vf),
\]
since $\vP\nabla_{\mathbf{w}}\vf^{\text{eq}}$ equals the identity matrix,
and then $\vtf$ satisfies

\[
\partial_{t}\vtf=\partial_{t}\vf-\nabla_{\mathbf{w}}\vf^{\text{eq}}(\vP\vf)\partial_{t}(\vP\vf)=0.
\]
This will lead to a specific time integration of the source term (see
Remark \ref{rem:source_G1}).

\end{rem}

\begin{rem}
System (\ref{eq:kin_bgk}) has to be supplemented with conditions
at the boundary $\partial\Omega$ of the computational domain $\Omega$.
We denote by $\normal=(n_{1}\ldots n_{\ddim})$ the outward normal
vector on $\partial\Omega.$ For simplicity, we shall only consider
very simple time-independent Dirichlet boundary conditions $\vf^{b}$.
We note
\[
\vV\cdot\normal=\sum_{k=1}^{\ddim}\vV^{k}n_{k},\quad\vV\cdot\normal^{+}=\max(\vV\cdot\normal,0),\quad\vV\cdot\normal^{-}=\min(\vV\cdot\normal,0).
\]
A natural boundary condition, which is compatible with the transport
operator, is 
\begin{equation}
\vV\cdot\normal^{-}\vf(\vx,t)=\vV\cdot\normal^{-}\vf^{b}(\vx),\quad\vx\in\partial\Omega.\label{eq:boundary_conditions-1}
\end{equation}

Boundary conditions (\ref{eq:boundary_conditions-1}) are very natural
from the kinetic point of view. However, they are not necessarily
natural when we go back to the macroscopic hyperbolic system. For
instance, at a given point of the boundary, the number of conditions
depends on the lattice velocities, which have no physical meaning.
It should rather depend on the number of characteristics of the macroscopic
system that are entering the computational domain. Then, it is not surprising
that we can observe instabilities arising from the boundary if we
apply the boundary condition (\ref{eq:boundary_conditions-1}). In
one of the test cases proposed in Section \ref{subsec:Flow-cylinder},
we will show how we can design appropriate boundary conditions when
the macroscopic model requires non-slip boundary conditions.

\end{rem}

\subsection{Examples}

Devising a kinetic approximation consists in giving the discrete velocities
and the projection matrix $\vP$ such that there exists a equilibrium
function $\vf^{\text{eq}}$ satisfying the compatibility conditions (\ref{eq:Pfeq_w})-(\ref{eq:flux_kin}).
We first present a generic method, the so-called vectorial kinetic
method, and then some specific Lattice-Boltzmann schemes.

\subsubsection{Vectorial kinetic method \label{subsec:vectorial-kinetic-scheme} }

The principle of the vectorial kinetic representation is to apply
an analogue decomposition to each component of the hyperbolic system
\cite{graille2014approximation}. 

We here present the simplest method belonging to this family. It consists
in choosing, for each component $w_{l}$ of macroscopic field $\mathbf{w}=(w_{1},\ldots,w_{m})$,
the same velocity set aligned with the Cartesian basis $(\mathbf{e}_{k},\ k=1,\dots D)$
and a unique velocity scale $\lambda$. For each  component $w_{l}$
of macroscopic field, we thus consider the $2D$ velocities 
\[
\mathbf{v}_{l,k,\pm}=\pm\lambda\mathbf{e}_{k},\quad k=1,\dots,D,
\]
and we note $f_{l,k,\pm}$ the corresponding components of the kinetic
distribution $\mathbf{f}$. We thus have $n_{v}=2D\times m$ discrete
velocities. 

The consistency conditions (\ref{eq:Pfeq_w}) and (\ref{eq:flux_kin})
yield $m\times(D+1)$ equations for the $2D\times m$ unknowns. The
projection $\mathbf{P}$ still remains to be defined. One possible
choice is to suppose that, for a given $l$ component, each $k$-th
velocity axis components ($f_{l,k,+}$ and $f_{l,k,-}$ ) contributes
to the macroscopic quantity $w_{l}$ in the same proportion. Hence,
relations (\ref{eq:Pfeq_w}) and (\ref{eq:flux_kin}) write
\begin{align*}
f_{l,k,+}^{\text{eq}}(\vw)+f_{l,k,-}^{\text{eq}}(\vw) & =\frac{w_{l}}{D},\quad\forall(l,k)\\
\lambda f_{l,k,+}^{\text{eq}}(\vw)-\lambda f_{l,k,-}^{\text{eq}}(\vw) & =q^{k}(\mathbf{w)}_{l},
\end{align*}

With this assumption, the equilibrium functions are uniquely defined
by

\begin{equation}
f_{l,j,\pm}^{\text{eq}}(\mathbf{w)}=\frac{w_{l}}{2D}\pm\frac{q^{j}(\mathbf{w)}_{l}}{2\lambda}.\label{eq:feq_cartesian}
\end{equation}

For these models, the diffusion tensor (\ref{eq:diffusion_vecto_general})
obtained from the Chapman-Enskog expansion takes a particular simple
form. Indeed, the components of the first part of the tensor simplify
into 

\begin{equation}
\left(P\vV^{k}\vV^{j}\nabla_{\mathbf{w}}\vf^{\text{eq}}\right)_{l,l'}=\frac{\lambda^{2}}{D}\delta_{kj}\delta_{l,l'},
\end{equation}

so that each directional block of the diffusion tensor writes

\begin{equation}
\mathcal{D}^{kj}=\frac{\lambda^{2}}{D}\text{\text{Id}}-\nabla_{\mathbf{w}}\vq^{k}\nabla_{\mathbf{w}}\vq^{j}.\label{eq:diffusion_vecto_cartesian}
\end{equation}

Then, from equation (\ref{eq:entropy}), considering a convex entropy
$\eta(\vw)$ of the macroscopic system, the limit system is entropy
dissipative provided the tensor

\begin{equation}
\sigma_{kj}=\nabla_{\vw}^{2}\eta\left[\frac{\lambda^{2}}{D}\text{Id}-\nabla_{\mathbf{w}}\vq^{k}\nabla_{\mathbf{w}}\vq^{j}\right]
\end{equation}

is definite non-negative. 
\begin{example}
\label{exa:Euler_1d}(One-dimensional isothermal Euler equations,
vectorial method) Let apply the above framework to the one-dimensional
isothermal compressible Euler equations. The conservative system is
given by $m=2$ and 
\begin{equation}
\v w=(\rho,\rho u)^{T},\label{eq:euler_w}
\end{equation}
\begin{equation}
\v q^{1}(\v w)=\v q(\v w)=(\rho u,\rho u^{2}+c^{2}\rho)^{T}.\label{eq:euler_q}
\end{equation}

where $\rho(x,t)$ is the density, $u(x,t)$ the velocity, and $c>0$
the sound speed, which is a given parameter. The vectorial kinetic
model is given by $n_{v}=4$ and
\[
\v V^{1}=\text{diag}(-\lambda,\lambda,-\lambda,\lambda),\quad\v P=\left(\begin{array}{cccc}
1 & 1 & 0 & 0\\
0 & 0 & 1 & 1
\end{array}\right),
\]
\[
f_{k,\pm}^{\text{eq}}=\frac{w_{k}}{2}\pm\frac{q(\v w)_{k}}{2\lambda},\quad k=1,2.
\]
The diffusion tensor reads 
\end{example}

\begin{equation}
\mathcal{D}^{11}=\left[\begin{array}{cc}
\lambda^{2}-(c^{2}-u^{2}) & -2u\\
-2u\left(c^{2}-u^{2}\right) & \lambda^{2}-(c^{2}+3u^{2})
\end{array}\right].
\end{equation}

An entropy for this system is $\eta=\rho\frac{u^{2}}{2}+\rho c^{2}\log(\frac{\rho}{\rho_{0}})$,
with $\rho_{0}>0$ an arbitrary constant. The entropy dissipation
tensor reads

\begin{equation}
\nabla_{\vw\vw}^{2}\eta\mathcal{D}^{11}=\frac{1}{\rho}\left[\begin{array}{cc}
2u^{2}(c^{2}-u^{2})+(c^{2}+u^{2})(\lambda^{2}-c^{2}+u^{2}) & -u(\lambda^{2}+c^{2}-u^{2})\\
-u(\lambda^{2}+c^{2}-u^{2}) & \lambda^{2}-c^{2}-u^{2}
\end{array}\right],
\end{equation}

which is definite non-negative provided $\lambda>\vert u\vert+c$.
The lattice velocity $\lambda$ has to satisfy the sub-characteristic
condition $\lambda>\left|u\right|+c.$ We note that this representation
is equivalent to the Jin and Xin relaxation \cite{jin1995relaxation}
of the associated hyperbolic system.

\subsubsection{Other Lattice Botzmann methods \label{subsec:LBM}}

We here present the D1Q3 Lattice Boltzmann scheme and the D2Q9 scheme,
its extension in two dimensions. 
\begin{example}\label{exa:ex4}
(One-dimensional isothermal Euler equations, D1Q3) The $D1Q3$ scheme
is a standard method for the one-dimensional isothermal Euler equations
(\ref{eq:euler_w})-(\ref{eq:euler_q}). This model takes advantage
of the structure of the Euler equations, which are moments of the
Boltzmann equation in the vanishing viscosity limit. The D1Q3 model
uses the velocity set $\mathbf{V}^{1}=(-\lambda,0,\lambda)$ and the
projection matrix 
\[
\vP=\left[\begin{array}{ccc}
1 & 1 & 1\\
-\lambda & 0 & \lambda
\end{array}\right].
\]
For this model, the diffusion tensor reads
\end{example}

\begin{equation}
\mathcal{D}^{11}(\vw)=\left[\begin{array}{cc}
0 & 0\\
-2u\left(c^{2}-u^{2}\right) & \lambda^{2}-3u^{2}-c^{2}
\end{array}\right].\label{eq:diff_tensor_D1Q3}
\end{equation}

A notable fact is that there is no diffusion on the density. Let us
now consider the same entropy $\eta=\rho\frac{u^{2}}{2}+\rho c^{2}\log(\frac{\rho}{\rho_{0}})$
as for the vectorial scheme. The entropy dissipation tensor reads

\begin{equation}
\sigma_{11}(\vw)=\nabla_{\vw}^{2}\eta(\vw)D^{11}(\vw)=\frac{1}{\rho}\left[\begin{array}{cc}
u^{2}(c^{2}-u^{2}) & -u(\lambda^{2}-c^{2}-3u^{2})\\
-u(c^{2}-u^{2}) & \lambda^{2}-c^{2}-3u^{2}
\end{array}\right].
\end{equation}

\begin{figure}
\begin{centering}
\includegraphics[width=0.3\textwidth]{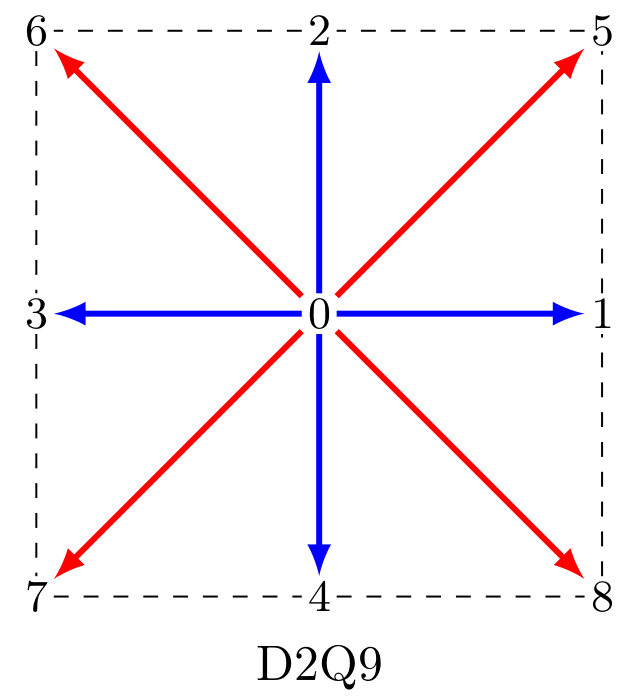}
\par\end{centering}
\caption{\label{fig:D2Q9_D3Q27_nodes}D2Q9  velocity grid.}
\end{figure}

Unfortunately neither $\mathcal{D}^{11}$ nor $\sigma_{11}$ can be
made definite positive by setting the value of $\lambda$. Indeed,
the symmetric part of $\sigma_{11}$ (resp. $\mathcal{D}^{11}$) has
always two real eigenvalues of opposite sign, regardless of the value
of $\lambda$. 
\begin{example}
\label{exa:2d_D2Q9}(Two-dimensional isothermal Euler equations, D2Q9)
The extension in two dimension of the previous D1Q3 scheme is the
D2Q9 scheme. The number of conservative variables for the two-dimensional
isothermal Euler scheme is $m=3$. The conservative variables are
\[
\v w=(\rho,\rho u,\rho v)^{T},
\]
and the flux is given
\[
\v q^{1}(\v w)=(\rho u,\rho u^{2}+c^{2}\rho,\rho uv)^{T},
\]
\[
\v q^{2}(\v w)=(\rho v,\rho uv,\rho v^{2}+c^{2}\rho)^{T},
\]
where the constant $c>0$ is the sound speed. The number of kinetic
equations is $n=9$. The kinetic model is based on a lattice of $n_{v}$
velocities $\boldsymbol{\v v}_{i}=(v_{i}^{1},v_{i}^{2})$, $i=1,\ldots,n_{v}$,
given by
\[
(v_{j}^{k})=\lambda\left(\begin{array}{ccccccccc}
0 & 1 & 0 & -1 & 0 & 1 & -1 & -1 & 1\\
0 & 0 & 1 & 0 & -1 & 1 & 1 & -1 & -1
\end{array}\right),
\]
with $\lambda>0$, and represented in Figure \ref{fig:D2Q9_D3Q27_nodes}.
The projection matrix is given by
\[
\v P=\left(\begin{array}{ccc}
1 & \cdots & 1\\
v_{1}^{1} & \cdots & v_{9}^{1}\\
v_{1}^{2} & \cdots & v_{9}^{2}
\end{array}\right).
\]
and the equilibrium distribution is given by
\[
\forall1\leqslant i\leqslant9,\text{\quad}f_{i}^{\text{eq}}(\vw)=\omega_{i}\,\rho\left(1+\frac{v_{i}^{1}u+v_{i}^{2}v}{c^{2}}+\frac{\left(v_{i}^{1}u+v_{i}^{2}v\right)^{2}}{c^{4}}-\frac{u^{2}+v^{2}}{2c^{2}}\right),
\]

where the weights are given by:

\[
\boldsymbol{\omega}=\left(\begin{array}{ccccccccc}
\frac{4}{9} & \frac{1}{9} & \frac{1}{9} & \frac{1}{9} & \frac{1}{9} & \frac{1}{36} & \frac{1}{36} & \frac{1}{36} & \frac{1}{36}\end{array}\right).
\]
In the literature, a common choice is to take $\lambda=\sqrt{3}c$.
In this case the kinetic model is stable only for low Mach number
flows (fluid velocity small compared to the sound speed).  It has nevertheless good properties (no diffusion on the density,
for instance) and a requires a small number of velocities. 
\end{example}

\section{Implicit high-order palindromic time-stepping\label{sec:Palindromic-time-stepping}}

In this section, we present a scheme for the kinetic-relaxation representation
(\ref{eq:kin_bgk}), which is implicit and high order in time. We
rewrite equation (\ref{eq:kin_bgk}) as follows

\begin{equation}
\partial_{t}\v f+\v L\v f+\v{\vN}\v f=0,\label{eq:kin_bgk_2}
\end{equation}

where $\vL\vf=\sum_{k=1}^{D}\vV^{k}\partial_{k}\vf$ , and $\vN\vf=-(\v f^{\text{eq}}(\vP\vf)-\v f)/\tau$.
The transport operator $\vL$ is linear, while the relaxation operator
$\vN$ is non-linear but local. For keeping the explanations simple,
we do not include the source term in the general presentation but
we will add some specific remarks.

For numerical applications, we shall consider an approximation $\v f_{h}$
of $\v f$ in a finite-dimensional space $E_{h}$. The parameter $h$
is for instance the size $\Delta x$ of the cells in the Discontinuous
Galerkin mesh. We assume that the approximation error behaves like
$O(h^{p})$ with $p\geq1$: the space approximation is at least first
order accurate with respect to the discretization parameter $h$.
The kinetic equation (\ref{eq:kin_bgk}) is thus approximated by a
set of differential equations
\begin{equation}
\partial_{t}\v f_{h}+\v L_{h}\v f_{h}+\vNh\v f_{h}=0,\label{eq:kinetic_approx}
\end{equation}
where the operator $\vNh$ actually equals the relaxation operator
$\vN$, but the operator $\v L_{h}$ is an approximation of $\vL$.
For the approximate transport operator $\vL_{h},$ several possibilities
may be considered: finite differences, finite elements, discrete Fourier
transform, Discontinous Galerkin (DG) approximation, semi-Lagrangian
methods, \textit{etc}. In this paper, we adopt an upwind nodal DG
approximation \cite{hesthaven2007nodal} (see Section \ref{sec:dG}). 

\subsection{First order splitting}

The exact flow of the differential equation (\ref{eq:kinetic_approx})
is given by
\[
\v f_{h}(t)=\exp(-t(\v L_{h}+\vNh))\,\v f_{h}(0).
\]

The exponential notation can be made completely rigorous here even
in the case of non-linear operators thanks to the Lie algebra formalism.
For an exposition of this formalism in the context of numerical methods
for ordinary differential equations, we refer for instance to \cite{hairer2006geometric,mclachlan2002splitting}.

Computing the exact flow is generally not possible. Instead, we apply
a splitting method in order to integrate the differential equation
(\ref{eq:kinetic_approx}). We can consider the simple Lie's splitting
approximation 

\begin{equation}
\vf_{h}(\Delta t)=\mathbf{M}_{1}(\Delta t)\vf_{h}(0)+O(\Delta t^{2}).\label{eq:liesplitting}
\end{equation}

with 
\[
\mathbf{M}_{1}(\Delta t)=\vR_{1}(\Delta t)\mathbf{T}_{1}(\Delta t)
\]

where $\vR_{1}$ and $\v T_{1}$ are first order approximations of
the relaxation and transport exact time integrators. In order to be
able to use large time steps, we consider the implicit first order
Euler scheme
\[
\vR_{1}(\Delta t)=(\text{{Id}}+\Delta t\,\vNh)^{-1},\quad\v T_{1}(\Delta t)=(\Id+\Delta t\,\v L_{h})^{-1}.
\]

For a fixed $\tau>0$, we actually have the estimates\footnote{For one single time step the error is $O(\Delta t^{2})$. But when
the error is accumulated on $t_{\max}/\Delta t$ time steps it indeed
produces a first order method.}
\[
\v R_{1}(\Delta t)=\exp(-\Delta t\,\vNh)+O(\Delta t^{2}),\quad\v T_{1}(\Delta t)=\exp(-\Delta t\v{\,L}_{h})+O(\Delta t^{2}).
\]
Let us point out that $\vR_{1}$ is a non-linear operator, because
$\v f\mapsto\v f^{\text{eq}}(\v P\vf)$ is non-linear. The linearity of $\v T_{1}$
depends on the linearity of $\v L_{h}$. The transport solver $\v L_{h}$
could be non-linear, even if the transport operator $\v{\vL}$ is
linear. This is the case if slope limiters are activated, for instance.

Finally, let us note that even if $\vR_{1}$ and $\v T_{1}$ are implicit
operators, they can actually be computed with an explicit cost. Indeed,
since for all $\vf$ we have $\v P\vNh\vf=0,$ the macroscopic quantity
$\mathbf{w=\vP}\vf$ is invariant during the relaxation step $\vR_{1}$.
It is then quite standard that $\vR_{1}$ takes the following explicit
form
\[
\vR_{1}(\Delta t)\v f=\frac{\v f^{\text{eq}}(\v{\vP f})+\frac{\tau}{\Delta t}\v f}{1+\frac{\tau}{\Delta t}}.
\]
In addition, because the free transport step is solved by an upwind
DG solver, then the linear operator $\Id+\Delta t\,\v L_{h}$ is block-triangular
\cite{badwaik2017task} and its inverse $\v T_{1}$ can also be computed
explicitly. We detail the method in Section \ref{sec:transport_dg}.
\begin{rem}
\label{rem:source_G1}When a source term is present in the model,
we further compose $\vR_{1}$ and $\v T_{1}$ with the following local
operator
\end{rem}

\[
\vG_{1}(\Delta t)=(\text{{Id}}+\Delta t\,\vH)^{-1}.
\]

where $\vH\vf_{h}=\nabla_{\v w}\vf^{\text{eq}}(\vP\vf_{h})\,\vs(\vP\vf_{h})$
is the kinetic source operator. Unlike $\vR_{1}$ and $\v T_{1}$,
this operator is a priori truly non-linear. However, as noticed in
Remark \ref{rem:source}, operator $\vS_{h}$ acts only on the macroscopic
variables $\vP\vf_{h}$. Consequently, we have

\[
\vG_{1}(\Delta t)\vf_{h}=\vf^{\text{eq}}\big(\vS_{1}(\Delta t)\vP\vf_{h}\big)+\Big(\vf_{h}-\vf^{\text{\text{eq}}}(\vP\vf_{h})\Big),
\]
where $\vS_{1}$ is the implicit Euler scheme on the macroscopic variables

\[
\vS_{1}(\Delta t)=(\text{{Id}}+\Delta t\,\vs)^{-1}.
\]

\subsection{Second-order stiffly accurate splitting}

Using methods of geometric integration {[}\ref{eq:liesplitting}{]},
we now consider a second-order in time scheme, that keeps second-order
accuracy in the limit $\tau\to0$.

We consider the second-order Crank-Nicolson scheme for the transport
equation
\begin{equation}
\mathbf{T}_{2}(\Delta t)=(\v{\Id}+\frac{\Delta t}{2}\v L_{h})(\v{\Id}-\frac{\Delta t}{2}\v L_{h})^{-1}.\label{eq:transport_cn}
\end{equation}
as well as for the relaxation operator
\[
\vR_{2}(\Delta t)=(\v{\Id}+\frac{\Delta t}{2}\vNh)(\v{\Id}-\frac{\Delta t}{2}\vNh)^{-1}.
\]

These transport and relaxation operators can be solved with the cost
of an explicit scheme. Indeed, since the macroscopic variables $\vw=\v P\vf$
is unchanged during the relaxation step, the relaxation operator (like
$\vR_{1})$ is only apparently implicit. We actually have the explicit
formula:
\begin{equation}
\vR_{2}(\Delta t)\vfh=\frac{(2\tau-\Delta t)\vfh}{2\tau+\Delta t}+\frac{2\Delta t\v{\,f}^{\text{eq}}(\v{\vP}\vfh)}{2\tau+\Delta t}.\label{eq:collision_cn}
\end{equation}

As regards the transport step $\mathbf{T}_{2}$, it involves an explicit
and an implicit transport both over a time interval $\Delta t/2$.
Like in the first order splitting, the implicit transport solution
can be computed at the cost of an explicit solver (see Section \ref{sec:transport_dg}).

If $\tau>0$, we observe that the operators $\mathbf{T}_{2}$ and
$\vR_{2}$ are \emph{time-symmetric}: if we set $\mathbf{O}_{2}=\mathbf{T}_{2}$
, $\mathbf{O}_{2}=\vR_{2}$, or $\mathbf{O}_{2}=\mathbf{S}_{2}$ a first order (or more) approximation,
then $\mathbf{O}_{2}$ satisfies
\begin{equation}
\mathbf{O}_{2}(-\Delta t)=\mathbf{O}_{2}(\Delta t)^{-1},\quad\mathbf{O}_{2}(0)=\Id.\label{eq:prop_sym}
\end{equation}

This property implies that, since $\mathbf{O}_{2}$ is necessarily a second
order approximation of the exact integrator \cite{mclachlan2002splitting,hairer2006geometric}.
Let us now note that when $\tau=0$, the relaxation operator becomes
independent of the time step and writes
\begin{equation}
\vRz_{2}(\Delta t)\vfh=2\v f^{\text{eq}}(\vP\vfh)-\vfh,\label{eq:max_sym}
\end{equation}
and then $\vRz_{2}$ does not satisfy (\ref{eq:prop_sym}) anymore.
However, we note that, due to the conservation of the macroscopic
variables, it is an involution

\begin{equation}
\vRz_{2}(\Delta t)\vRz_{2}(\Delta t)=\Id.\label{eq:max_sym-1}
\end{equation}
This is the key point of the following scheme.

We propose to use the following time-symmetric splitting 

\begin{equation}
\vM_{2}(\Delta t)=\vT_{2}\left(\frac{\Delta t}{4}\right)\vR_{2}\left(\frac{\Delta t}{2}\right)\vT_{2}\left(\frac{\Delta t}{2}\right)\vR_{2}\left(\frac{\Delta t}{2}\right)\vT_{2}\left(\frac{\Delta t}{4}\right).
\end{equation}
It can be easily checked that $\vM_{2}(\Delta t)$ is time-symmetric
for all $\tau\geqslant0$ , including the case $\tau=0$. Consequently,
the scheme remains second order accurate in the limit $\tau\rightarrow0$.
\begin{rem}
The classical second-order Strang splitting,
\[
\tilde{\vM}_{2}(\Delta t)=\mathbf{T}_{2}\left(\frac{\Delta t}{2}\right)\vR_{2}(\Delta t)\mathbf{T}_{2}\left(\frac{\Delta t}{2}\right),
\]
is time-symmetric for all $\tau>\text{0}$ but not for $\tau=0$.
However for $\tau=0$, unless the method does not give the identity
operator on the kinetic distribution for $\Delta t=0$, it turns out
to be the identity operator on the macroscopic variables: $\vP\tilde{\vM}_{2}(0)\mathbf{\vf}=\vP\mathbf{\vf}$.
This might explain why second-order accuracy can be numerically observed
at $\tau=0$ for the macroscopic variables, even though the operator
on the full kinetic system is not symmetric. 
\end{rem}

\begin{rem}
To take into account source terms, we consider the following second-order
scheme
\[
\vG_{2}(\Delta t)\vf_{h}=\vf^{\text{eq}}\big(\vS_{2}(\Delta t)\vP\vf_{h}\big)+\Big(\vf_{h}-\vf^{\text{\text{eq}}}(\vP\vf_{h})\Big),
\]
where $\vS_{2}$ is the (truly) implicit Crank-Nicolson scheme on
the macroscopic variables
\[
\vS_{2}(\Delta t)=(\text{{Id}}+\Delta t\,\vs)(\text{{Id}}+\Delta t\,\vs)^{-1}.
\]
Then the second-order splitting is modified into
\begin{align*}
 & \bar{\vM}_{2}(\Delta t)=\\
 & \quad\vT_{2}\left(\frac{\Delta t}{4}\right)\vG_{2}\left(\frac{\Delta t}{2}\right)\vR_{2}\left(\frac{\Delta t}{2}\right)\vT_{2}\left(\frac{\Delta t}{2}\right)\vR_{2}\left(\frac{\Delta t}{2}\right)\vG_{2}\left(\frac{\Delta t}{2}\right)\vT_{2}\left(\frac{\Delta t}{4}\right).
\end{align*}

which is still time-symmetric.
\end{rem}

\subsection{\label{subsec:Palindromic-time-integration}High-order palindromic
splitting}

Once defined a second-order accurate time-symmetric scheme, palindromic
composition method enables to easily achieve any even order of accuracy
\cite{mclachlan2002splitting,hairer2006geometric,coulette2016palindromic}.
A general palindromic scheme with $s+1$ steps has the form
\begin{equation}
\vM{}_{p}(\Delta t)=\vM{}_{2}(\gamma_{0}\Delta t)\vM_{2}(\gamma_{1}\Delta t)\cdots\vM_{2}(\gamma_{s}\Delta t),\label{eq:palindromic}
\end{equation}
where the $\gamma_{i}$'s are real numbers satisfying
\[
\gamma_{i}=\gamma_{s-i},\quad0\leq i\leq s.
\]

In the following, we will consider the fourth-order Suzuki scheme
\cite{suzuki1990fractal,hairer2006geometric,mclachlan2002splitting}
and the sixth-order Kahan-Li scheme \cite{kahan1997composition},
whose intermediate steps are given in Table \ref{tab:Palindromic-coefficients}.
The Sukuki scheme requires $5$ steps, while the Kahan-Li scheme is
made of $9$ steps. 

We note that the two methods require to apply the elementary relaxation
or transport $\vR_{2}$ and $\mathbf{\vT}_{2}$ with negative time
steps $-\Delta t<0$. If we were using the exact transport solver
$\vL$, negative time steps would not cause any problem. However,
the transport approximation $\vL_{h}$ generally introduces a slight
dissipation to ensure stability (for instance upwinding in DG discretization,
see Section \ref{sec:dG}). In order to ensure stability, we have
thus to replace $\mathbf{T}_{2}(-\Delta t)$ with a more stable operator.
This can be done by observing that solving $\partial_{t}\vf+\vL\vf=0$
for negative time $t<0$ is equivalent to solve $\partial_{t'}\vf-\vL\vf=0$
for $t'=-t>0$. Therefore, we use 
\[
\mathbf{T}'_{2}(\Delta t)=(\v{\Id}+\frac{\Delta t}{2}(\v{-L})_{h})(\v{\Id}-\frac{\Delta t}{2}(-\v L)_{h})^{-1}
\]

where $(-\vL)_{h}$ is a stable discretization of $-\vL$. The numerical
relaxation operator $\vR_{2}$ is time reversible in the limit $\tau\to\text{0}$:
for $\tau=0$, it actually does not depend on $\Delta t$ anymore
(see (\ref{eq:max_sym})). In this stage, negative time steps do not
cause any difficulty, at least when $\tau\ll\Delta t$.

\begin{table}
$\begin{array}{c}
\text{\text{Suzuki coefficients} (p=4, s=4)}\\
\gamma_{0}=\gamma_{1}=\gamma_{3}=\gamma_{4}=\frac{1}{4-4^{1/3}},\quad\gamma_{2}=-\frac{4^{1/3}}{4-4^{1/3}}.\\
\text{\text{Kahan-Li coefficients} (p=6, s=8)}\\
\begin{array}{cc}
\gamma_{0}=\gamma_{8}= & 0.392161444007314139275655330038\ldots\\
\gamma_{1}=\gamma_{7}= & 0.332599136789359438604272125325\ldots\\
\gamma_{2}=\gamma_{6}= & -0.7062461725576393598098453372227\ldots\\
\gamma_{3}=\gamma_{5}= & 0.0822135962935508002304427053341\ldots\\
\gamma_{4}= & 0.798543990934829963398950353048\ldots
\end{array}
\end{array}$

\caption{Palindromic coefficients.\label{tab:Palindromic-coefficients}}
\end{table}

\section{Implicit discontinuous Galerkin method for linear transport\label{sec:dG}}

\label{sec:transport_dg} In this section, we briefly present the
approximate linear transport operator $\vL_{h}$ obtained from the
Discontinuous Galerkin method. We also show how its matrix-triangular
structure enables to solve, with an explicit cost, the implicit operator
involved in the second order Crank-Nicolson solver (see (\ref{eq:transport_cn})). 

\subsection{DG approximation}

As said above, for solving (\ref{eq:kin_bgk}) we treat the transport
operator $\vV\cdot\boldsymbol{\partial}$ and the collision operator
$\v N$ separately, thanks to the splitting approach. Let us now describe
the transport solver.

For a simple exposition, we only consider one single scalar transport
equation for $f(\vx,t)\in\mathbb{R}$ at constant velocity $\mathbf{v}$
\begin{equation}
\partial_{t}f+\mathbf{v}\cdot\nabla f=0.\label{eq:single_transport}
\end{equation}
The general vectorial case is easily deduced.

We consider a mesh $\mathcal{M}$ of $\Omega$ made of open sets,
called ``cells'', $\mathcal{M}=\left\{ L_{i},\,i=1\ldots N_{c}\right\} $.
In the most general setting, the cells satisfy
\begin{enumerate}
\item $L_{i}\cap L_{j}=\emptyset$, if $i\neq j$;
\item $\overline{\cup_{i}L_{i}}=\overline{\Omega}.$
\end{enumerate}
In each cell $L\in\mathcal{M}$, we consider a basis of functions $(\varphi_{L,i}(\vx))_{i=0\ldots N_{\dorder}-1}$
constructed from polynomials of order $\dorder$. We denote by $h$
the maximal diameter of the cells. With an abuse of notation we still
denote by $f$ the approximation of $f$, defined by
\[
f(\vx,t)=\sum_{j=0}^{N_{\dorder}-1}f_{L,j}(t)\varphi_{L,j}(\vx),\quad\vx\in L.
\]
The DG formulation then reads: find the $f_{L,j}$'s such that for
all cell $L$ and all test function $\varphi_{L,i}$
\begin{equation}
\int_{L}\partial_{t}f\varphi_{L,i}-\int_{L}f\mathbf{v}\cdot\nabla\varphi_{L,i}+\int_{\partial L}(\mathbf{v}\cdot\normal^{+}f_{L}+\mathbf{v}\cdot\normal^{-}f_{R})\varphi_{L,i}=0.\label{eq:dg_var}
\end{equation}
In this formula (see Figure \ref{eq:convention_downwind}):
\begin{itemize}
\item $R$ denotes the neighbouring cell to $L$ along its boundary $\partial L\cap\partial R$,
or the exterior of $\Omega$ on $\partial L\cap\partial\Omega$.
\item $\normal=\normal_{LR}$ is the unit normal vector on $\partial L$
oriented from $L$ to $R$.
\item $f_{R}$ denotes the value of $f$ in the neighbouring cell $R$ on
$\partial L\cap\partial R$.
\item If $L$ is a boundary cell, one may have to use the boundary values
instead: $f_{R}=f^{b}$ on $\partial L\cap\partial\Omega$.
\item $\mathbf{v}\cdot\normal^{+}f_{L}+\mathbf{v}\cdot\normal^{-}f_{R}$
is the standard upwind numerical flux encountered most finite volume
or DG methods.
\end{itemize}
\begin{figure}
\begin{centering}
\includegraphics[width=0.3\textwidth]{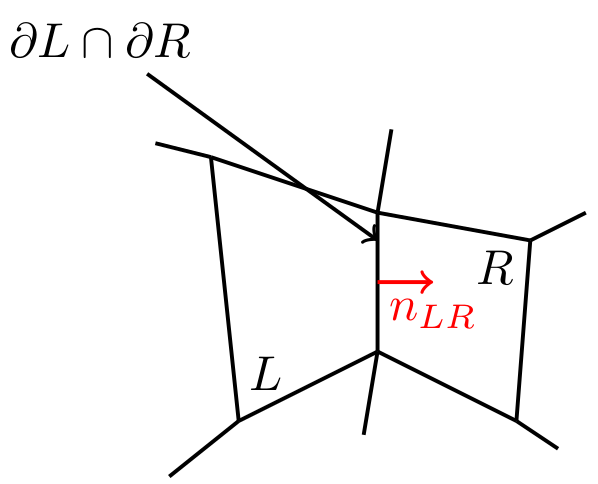}
\par\end{centering}
\caption{\label{fig:cell_convention}Convention for the $L$ and $R$ cells
orientation.}
\end{figure}

In our applications, we consider hexahedral cells. We have a reference
cell 
\[
\hat{L}=]-1,1[^{\ddim}
\]
\global\long\def\jacob{\boldsymbol{\tau}}
and a smooth transformation $\vx=\jacob_{L}(\hat{\vx})$, $\hat{\vx}\in\hat{L}$,
that maps $\hat{L}$ on $L$
\[
\jacob_{L}(\hat{L})=L.
\]
We assume that $\jacob_{L}$ is invertible and we denote by $\jacob_{L}'$
its (invertible) Jacobian matrix. We also assume that $\jacob_{L}$
is a direct transformation\textcompwordmark{}
\[
\det\jacob_{L}'>0.
\]

In our implementation $\jacob_{L}$ is a quadratic map based on hexahedral
curved ``H20'' finite elements with 20 nodes. The mesh of H20 finite
elements is generated by \texttt{gmsh} \cite{geuzaine2009gmsh}. 

On the reference cell, we consider the Gauss-Lobatto (GL) points $(\hat{\vx}_{i})_{i=0\ldots N_{\dorder}-1}$,
$N_{\dorder}=(\dorder+1)^{\ddim}$ and associated weights $(\omega_{i})_{i=0\ldots N_{\dorder-1}}$.
They are obtained by tensor products of the $(\dorder+1)$ one-dimensional
Gauss-Lobatto (GL) points on $]-1,1[$. The reference GL points and
weights are then mapped to the physical GL points of cell $L$ by
\begin{equation}
\vx_{L,i}=\jacob_{L}(\hat{\vx}_{i}),\quad\omega_{L,i}=\omega_{i}\det\jacob_{L}'(\hat{\vx}_{i})>0.\label{eq:map_GL}
\end{equation}
In addition, the six faces of the reference hexahedral cell are denoted
by $F_{\epsilon}$, $\epsilon=1\ldots6$ and the corresponding outward
normal vectors are denoted by $\hat{\normal}_{\epsilon}$. A big advantage
of choosing the GL points is that the volume and the faces share the
same quadrature points. A special attention is necessary for defining
the face quadrature weights. If a GL point $\hat{\vx}_{i}\in F_{\epsilon}$,
we denote by $\mu_{i}^{\epsilon}$ the corresponding quadrature weight
on face $F_{\epsilon}$. We also use the convention that $\mu_{i}^{\epsilon}=0$
if $\hat{\vx}_{i}$ does not belong to face $F_{\epsilon}$. A given
GL point $\hat{\vx}_{i}$ can belong to several faces when it is on
an edge or in a corner of $\hat{L}$. Because of symmetry, we observe
that if $\mu_{i}^{\epsilon}\neq0$, then the weight $\mu_{i}^{\epsilon}$
does not depend on $\epsilon$.

We then consider basis functions $\hat{\varphi_{i}}$ on the reference
cell: they are the Lagrange polynomials associated to the Gauss-Lobatto
point and thus satisfy the interpolation property
\[
\hat{\varphi}_{i}(\hat{\vx}_{j})=\delta_{ij}.
\]
The basis functions on cell $L$ are then defined according to the
formula
\[
\varphi_{L,i}(\vx)=\hat{\varphi}_{i}(\jacob_{L}^{-1}(\vx)).
\]
In this way, they also satisfy the interpolation property
\begin{equation}
\varphi_{L,i}(\vx_{L,j})=\delta_{ij}.\label{eq:interp_property}
\end{equation}
In this paper, we only consider conformal meshes: the GL points on
cell $L$ are supposed to match the GL points of cell $R$ on their
common face.

Let $L$ and $R$ be two neighbouring cells. Let $\vx_{L,j}$ be a
GL point in cell $L$ that is also on the common face between $L$
and $R$. In the case of conformal meshes, it is possible to define
the index $j'$ such that
\[
\vx_{L,j}=\vx_{R,j'}.
\]

Applying a numerical integration to (\ref{eq:dg_var}), using (\ref{eq:map_GL})
and the interpolation property (\ref{eq:interp_property}), we finally
obtain
\begin{multline}
\partial_{t}f_{L,i}\omega_{L,i}-\sum_{j=0}^{N_{\dorder}-1}\mathbf{v}\cdot\nabla\varphi_{L,i}(\vx_{L,j})f_{L,j}\omega_{L,j}+\\
\sum_{\epsilon=1}^{6}\mu_{i}^{\epsilon}\left(\mathbf{v}\cdot\normal_{\epsilon}(\vx_{L,i})^{+}f_{L,i}+\mathbf{v}\cdot\normal_{\epsilon}(\vx_{L,i})^{-}f_{R,i'}\right)=0.\label{eq:dg_lobatto}
\end{multline}
We have to detail how the gradients and normal vectors are computed
in the above formula. Let $\v A$ be a square matrix. We recall that
the cofactor matrix of $\v A$ is defined by
\begin{equation}
\text{co}(\v A)=\det(\v A)\left(\v A^{-1}\right)^{T}.\label{eq:co_mat}
\end{equation}
The gradient of the basis function is computed from the gradients
on the reference cell using (\ref{eq:co_mat})
\[
\nabla\varphi_{L,i}(\vx_{L,j})=\frac{1}{\det\jacob_{L}'(\hat{\vx}_{i})}\text{co}(\jacob_{L}'(\hat{\vx}_{j}))\hat{\nabla}\hat{\varphi}_{i}(\hat{\vx}_{j}).
\]
In the same way, the scaled normal vectors $\v n_{\epsilon}$ on the
faces are computed by the formula
\[
\normal_{\epsilon}(\vx_{L,i})=\text{co}(\jacob_{L}'(\hat{\vx}_{i}))\hat{\normal}_{\epsilon}.
\]
We introduce the following notation for the cofactor matrix\global\long\def\comat{\mathbf{c}}
\[
\comat_{L,i}=\text{co}(\jacob_{L}'(\hat{\vx}_{i})).
\]
The nodal DG scheme then reads
\begin{multline}
\partial_{t}f_{L,i}-\frac{1}{\omega_{L,i}}\sum_{j=0}^{N_{\dorder}-1}\mathbf{v}\cdot\comat_{L,j}\hat{\nabla}\hat{\varphi}_{i}(\hat{\vx}_{j})f_{L,j}\omega_{j}+\\
\frac{1}{\omega_{L,i}}\sum_{\epsilon=1}^{6}\mu_{i}^{\epsilon}\left(\mathbf{v}\cdot\comat_{L,i}\hat{\normal}_{\epsilon}{}^{+}f_{L,i}+\mathbf{v}\cdot\comat_{L,i}\hat{\normal}_{\epsilon}{}^{-}f_{R,i'}\right)=0.\label{eq:DG_reduit}
\end{multline}
On boundary GL points, the value of $f_{R,i'}$ is given by the boundary
condition
\[
f_{R,i'}=f^{b}(\vx_{L,i}),\quad\vx_{L,i}=\vx_{R,i'}.
\]
For practical reasons, it is interesting to also consider $f_{R,i'}$
as an artificial unknown in the fictitious cell. The fictitious unknown
is then a solution of the differential equation
\begin{equation}
\partial_{t}f_{R,i'}=0.\label{eq:fictitious_ODE}
\end{equation}
In the end, if we put all the unknowns in a single vector $\v F(t)$,
(\ref{eq:DG_reduit}), (\ref{eq:fictitious_ODE}) read as a large
system of coupled differential equations\global\long\def\transmat{\mathbf{L}_{h}}
\begin{equation}
\partial_{t}\v F_{h}=\transmat\v F_{h}.\label{eq:linear_ODE}
\end{equation}

This defines $\transmat$ the transport matrix. The transport matrix
satisfies the following properties:
\begin{itemize}
\item $\transmat\v F_{h}=0$ if the components of $\v F$ are all the same.
\item Let $\v F_{h}$ be such that the components corresponding to the boundary
terms vanish. Then for the scalar product $\left\langle \v F,\v G\right\rangle =\sum_{L}\sum_{i}\omega_{L,i}f_{L,i}g_{L,i}\,,$
we have 
\begin{equation}
\left\langle \v F_{h},\transmat\v F_{h}\right\rangle \leq0.\label{eq:dissipation}
\end{equation}
This dissipation property is a consequence of the choice of an upwind
numerical flux \cite{johnson1984finite}\footnote{Actually, this dissipation property is true only when the geometrical
transformations $\jacob_{L}$ are affine maps. For quadratic maps,
the Gauss-Lobatto numerical integration is not exact anymore (``aliasing''
effect: see \cite{hesthaven2007nodal} for instance). Weak instabilities
may develop for long-time numerical simulations.}.
\item In many cases, and with a good numbering of the unknowns in $\v F_{h}$,
$\transmat$ has a block-triangular structure. This aspect is discussed
in Subsection \ref{subsec:triangular_subsection}.
\end{itemize}
As stated above, we actually have to apply a transport solver for
each constant velocity $\v v_{i}$.

Let $L$ be a cell of the mesh $\mathcal{M}$ and $\vx_{i}$ a GL
point in $L$. As in the scalar case, we denote by $\vf_{L,i}$ the
approximation of $\vf$ in $L$ at GL point $i$. In the sequel, with
an abuse of notation and according to the context, we may continue
to note $\v F(t)$ the big vector made of all the vectorial values
$\vf_{L,j}$ at all the GL points $j$ in all the (real or fictitious)
cells $L$.

We may also continue to denote by $\v L_{h}$ the matrix made of the
assembly of all the transport operators for all velocities $\v v_{i}$.
With a good numbering of the unknowns it is possible in many cases
to suppose that $\v L_{h}$ is block-triangular. More precisely, because
in the transport step the equations are uncoupled, we see that $\v L_{h}$
can be made block-diagonal, each diagonal block being itself block-triangular.
See next Section \ref{subsec:triangular_subsection}.

\subsection{\label{subsec:triangular_subsection}Triangular structure of the
transport matrix}

\begin{figure}
\includegraphics[height=6cm]{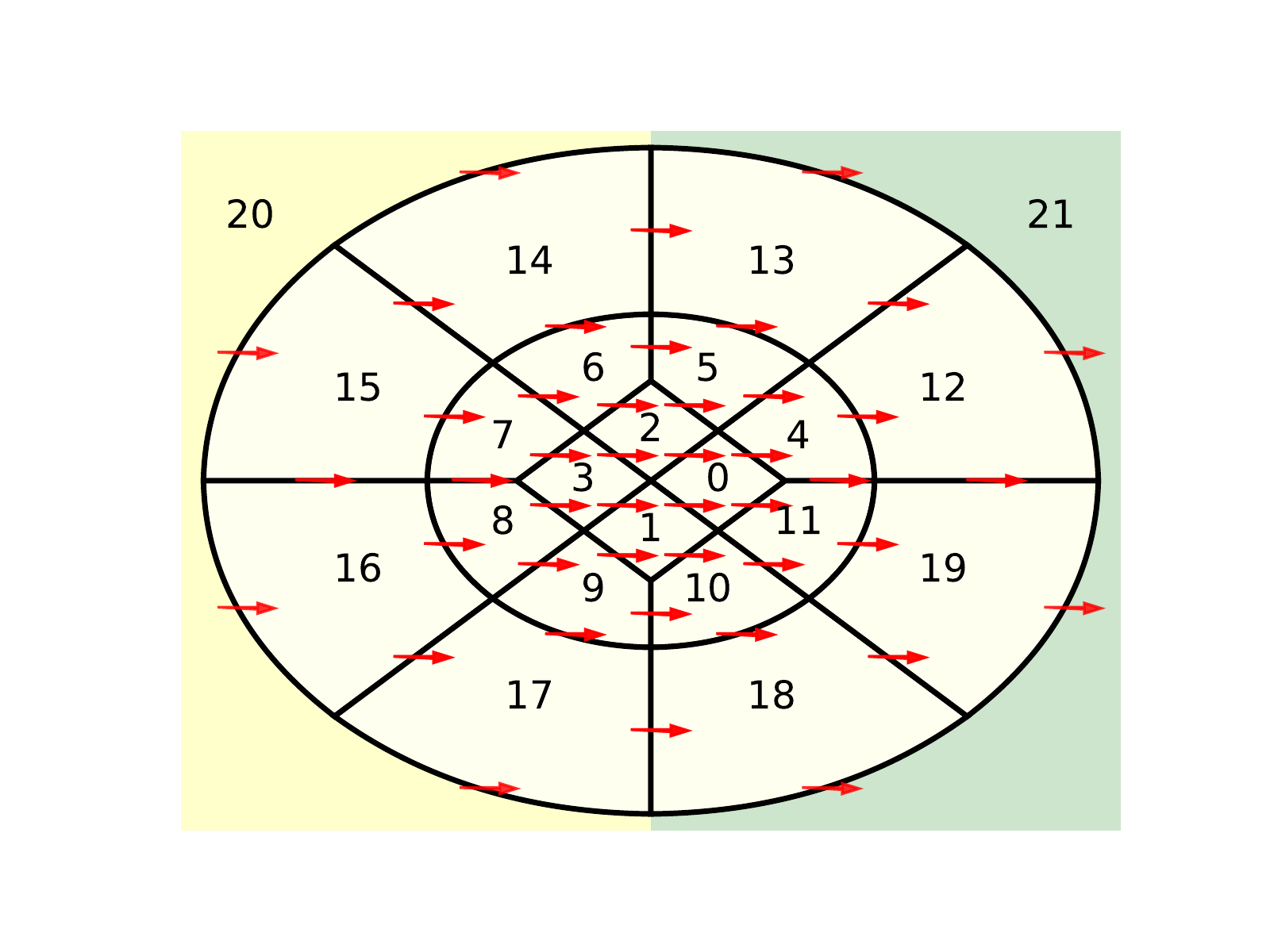}\includegraphics[height=10cm]{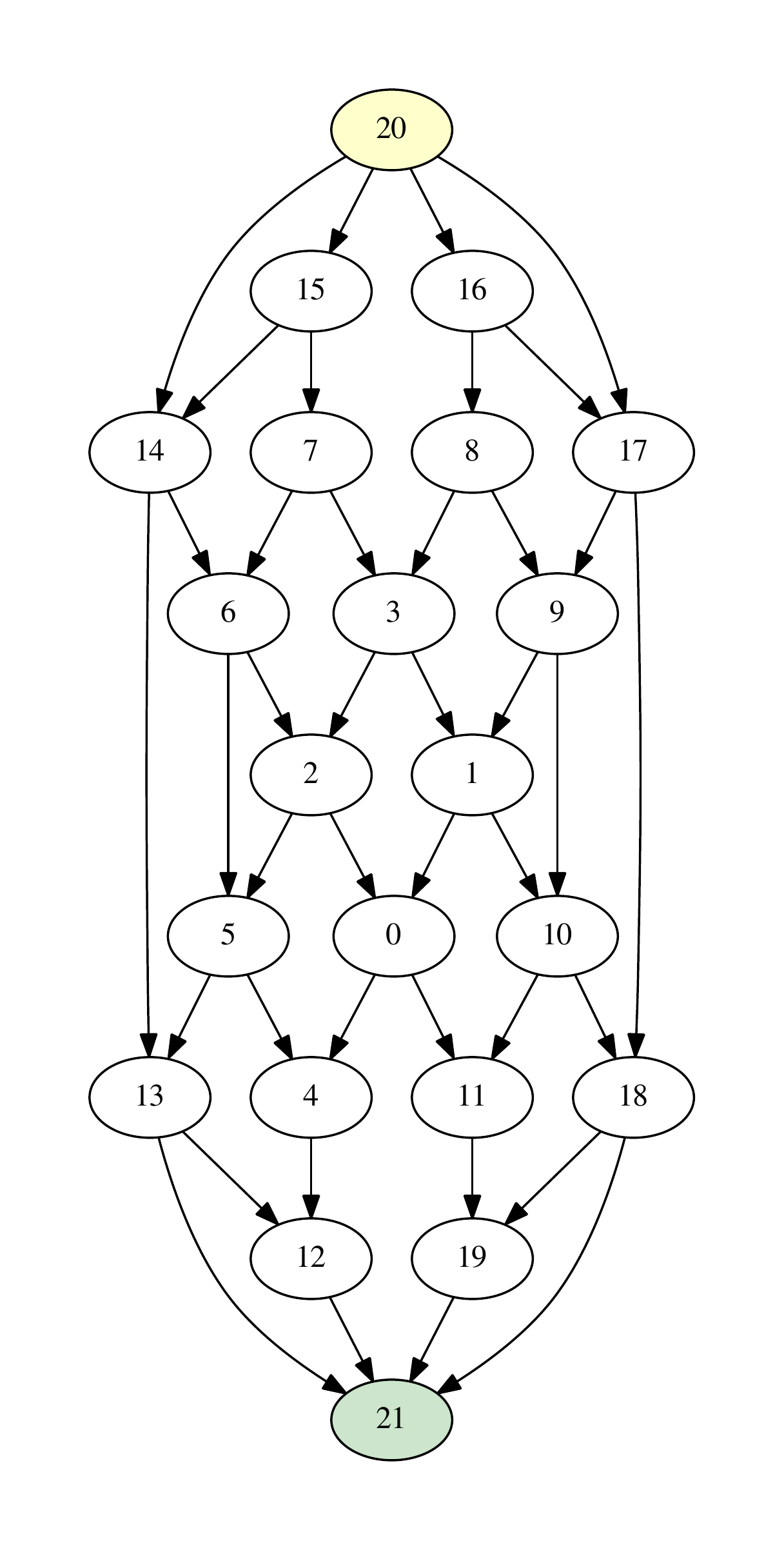}

\caption{\label{fig:dependency_graph}Construction of the dependency graph.
Left: example of a simple unstructured mesh with 20 interior cells.
The velocity field $v$ is indicated by red arrows. We add two fictitious
cells: one for the upwind boundary condition (cell $20$) and one
for the outflow part of $\partial\Omega$ (cell $21$). Right: the
corresponding dependency graph $\mathcal{G}.$ By examining the dependency
graph, we observe that the values of $\protect\v F^{n+1}$ in cell
$15$ and $16$ have to be computed first, using the boundary conditions.
Then cells $[14,7,8,17]$ can be computed in parallel, then cells
$[6,3,9]$ can be computed in parallel, and so on until the downwind
boundary is reached.}
\end{figure}

Because of the upwind structure of the numerical flux, it appears
that the transport matrix is often block-triangular. This is very
interesting because this allows to applying implicit schemes to (\ref{eq:linear_ODE})
without the costly inversion of linear systems \cite{moustafa20153d}.
We can provide the formal structure of $\transmat$ through the construction
of a directed graph $\mathcal{G}$ with a set of vertices $\mathcal{V}$
and a set of edges $\mathcal{E}\subset\mathcal{V}\times\mathcal{V}$.
The vertices of the graph are associated to the (real or fictitious)
cells of $\mathcal{M}$. Now consider two cells $L$ and $R$ with
a common face $F_{LR}$. We denote by $\normal_{LR}$ the normal vector
on $F_{LR}$ oriented from $L$ to $R$. If there is at least one
GL point $\vx$ on $F_{LR}$ such that 
\[
\normal_{LR}(\vx)\cdot\mathbf{v}>0,
\]
then the edge from $L$ to $R$ belongs to the graph:
\[
(L,R)\in\mathcal{E},
\]
see Figure \ref{fig:dependency_graph}.

In (\ref{eq:DG_reduit}) we can distinguish between several kinds
of terms. We write
\[
\partial_{t}f_{L}+\Gamma_{L\leftarrow L}f_{L}+\sum_{(R,L)\in\mathcal{E}}\Gamma_{L\leftarrow R}f_{R},
\]
with
\begin{multline*}
\Gamma_{L\leftarrow L}f_{L}=-\frac{1}{\omega_{L,i}}\sum_{j=0}^{N_{\dorder}-1}\mathbf{v}\cdot\comat_{L,j}\hat{\nabla}\hat{\varphi}_{i}(\hat{\vx}_{j})f_{L,j}\omega_{j}+
\frac{1}{\omega_{L,i}}\sum_{\epsilon=1}^{6}\mu_{i}^{\epsilon}\mathbf{v}\cdot\comat_{L,i}\hat{\normal}_{\epsilon}{}^{+}f_{L,i},
\end{multline*}
and, if $(R,L)\in\mathcal{E}$,
\[
\Gamma_{L\leftarrow R}f_{R}=\frac{1}{\omega_{L,i}}\mu_{i}^{\epsilon}\mathbf{v}\cdot\comat_{L,i}\hat{\normal}_{\epsilon}{}^{-}f_{R,i'}.
\]
We can use the following convention
\begin{equation}
(R,L)\notin\mathcal{E}\Rightarrow\Gamma_{L\leftarrow R}=0.\label{eq:convention_downwind}
\end{equation}
$\Gamma_{L\leftarrow L}$ contains the terms that couple the values
of $f$ inside the cell $L$. They correspond to diagonal blocks of
size $(\dorder+1)^{\ddim}\times(\dorder+1)^{\ddim}$ in the transport
matrix $\transmat$. $\Gamma_{L\leftarrow R}$ contains the terms
that couple the values inside cell $L$ with the values in the neighboring
upwind cell $R$. If $R$ is a downwind cell relatively to $L$ then
$\mu_{i}^{\epsilon}\mathbf{v}\cdot C_{L,i}\hat{\normal}_{\epsilon}{}^{-}=0$
and $\Gamma_{L\leftarrow R}=0$ is indeed compatible with the above
convention (\ref{eq:convention_downwind}).

Once the graph $\mathcal{G}$ is constructed, we can analyze it with
standard tools. If it contains no cycle, then it is called a Directed
Acyclic Graph (DAG). Any DAG admits a topological ordering of its
nodes. A topological ordering is a numbering of the cells $i\mapsto L_{i}$
such that if there is a path from $L_{i}$ to $L_{j}$ in $\mathcal{G}$
then $j>i$. In practice, it is useful to remove the fictitious cells
from the topological ordering. In our implementation they are put
at the end of the list.

Once the new ordering of the graph vertices is constructed, we can
construct a numbering of the components of $\v F$ by first numbering
the unknowns in $L_{0}$ then the unknowns in $L_{1}$, \emph{etc}.
More precisely, we set
\[
F_{kN_{\dorder}+i}=f_{L_{k},i}.
\]
Then, with this ordering, the matrix $\transmat$ is lower block-triangular
with diagonal blocks of size $(\dorder+1)^{\ddim}\times(\dorder+1)^{\ddim}$.
It means that we can apply implicit schemes to (\ref{eq:linear_ODE})
without costly inversion of large linear systems.

As stated above, we actually have to apply a transport solver for
each constant velocity $\v v_{i}$. In the sequel, with another abuse
of notation and according to the context, we continue to note $\v F$
the big vector made of all the vectorial values $\vf_{L,j}$ at all
the GL points $j$ in all the (real or fictitious) cells $L$.

We may also continue to denote by $\v L_{h}$ the matrix made of the
assembly of all the transport operators for all velocities $\v v_{i}$.
With a good numbering of the unknown it is still possible to suppose
that $\v L_{h}$ is block-triangular. More precisely, as in the transport
step the equations are uncoupled, we see that $\v L_{h}$ can be made
a block-diagonal matrix, each diagonal block being itself block-triangular.

\section{Parallel implementation}

Thanks to the splitting procedure described in Section \ref{sec:Palindromic-time-stepping}
the whole algorithm exhibits several levels of parallelism. First,
it is clear that the collision step is purely local to each interpolation
point and thus embarrassingly parallel. Second, the transport equations
are completely uncoupled from the other ones. They can thus be solved
independently in parallel. Finally, as stated above (see Figure \ref{fig:dependency_graph}),
inside the resolution of each transport equation it is again possible
to detect additional parallelism from the examination of the dependency
graph.

We have written a C implementation of the Palindromic Discontinuous Galerkin (PDG) method using a data-based
formulation of the parallelism. In this formulation it is essential
to distinguish between the input (Read mode) and output (Write mode)
data of each elementary computational task. The tasks are then submitted
to a runtime system that is able to distribute the work on the available
processors. From the data dependency, the runtime system detects the
tasks that can be performed in parallel. In our implementation, we
rely on the StarPU runtime library, which is especially designed for
efficient scientific computing \cite{AugAumFurNamThi2012EuroMPI}.
We use the MPI version of StarPU in order to distribute the computations
on clusters of multicore computers. 

Submitting a task to the StarPU system induces a slight overhead.
It is thus important to submit tasks that are not too small (too much
time would be spent into the tasks management) or not too big (which
could block the  tasks flow). Therefore, we apply what we call a ``macrocell''
approach. The geometry is first meshed at a coarse level. We call
the cells of the coarse mesh the ``macrocells''. The macrocells
are then refined into several subcells. We apply the task-based transport
solver described in Figure \ref{fig:dependency_graph} at the macrocell
level instead of the subcell level. In this way, we can adjust the
grain of the parallelism. This approach necessitates solving local
transport equations into the macrocells. This is achieved by assembling
and solving local block-triangular linear system. Those local sparse
linear systems are solved with the KLU library, which is able to detect
efficiently block-triangular structures \cite{davis2010algorithm}.
More details on the implementation are given in \cite{badwaik2017task}.

For the moment, the local systems are assembled and factorized at
each time-step. It would probably be more efficient to store the local
LU decompositions for saving computational time. We have not yet compared
the efficiency of our approach with other explicit or implicit DG
solvers. However, we have observed a good parallel scaling of the
method when the number of computational cores increases \cite{badwaik2017task}.
In addition, as it is shown in the numerical sections, the PDG method
accepts very high CFL numbers, which makes it a good candidate for
avoiding costly non-linear implicit solvers.

\section{Numerical results}

In this section, we apply the methodology presented in the previous
sections. We first numerically demonstrate the accuracy of the scheme
on one-dimensional test cases. We then show how the method applies
to two-dimensional models. We will make some remarks on the treatment
of the boundary conditions.

An important feature of the PDG method is the possibility to consider
large time steps without oscillations. In order to measure this advantage,
we have to define precisely how we define the time step and the corresponding
CFL number.

\subsection{One-dimensional isothermal Euler test cases}

In this section, we consider the vectorial kinetic method apply to
the one-dimensional isothermal Euler system, presented in Example
\ref{exa:Euler_1d}.

\subsubsection{Smooth solution}

For the first validation of the method we consider a test case with
a smooth solution, in the fluid limit $\tau=0$. The initial condition
is given by
\[
\rho(x,0)=1+e^{-30x^{2}},\quad u(x,0)=0.
\]
The sound speed is set to $c=0.6$ and the lattice velocity to $\lambda=2$.
We define the CFL number $\beta=\lambda\Delta t/\delta$, where $\delta$
is the minimal distance between two Gauss-Lobatto points in the mesh.
First, the CFL number is fixed to $\beta=5$. We consider a sufficiently
large computational domain $[a,b]=[-2,2]$ and a sufficiently short
final time $t_{\max}=0.4$ so that the boundary conditions play no
role. The reference solution $\vf(\cdot,t_{\max})$ is computed numerically
with a very fine mesh. In the DG solver the polynomial order in $x$
is fixed to $d=5$.

\begin{figure}
\centering{}\includegraphics[width=6cm]{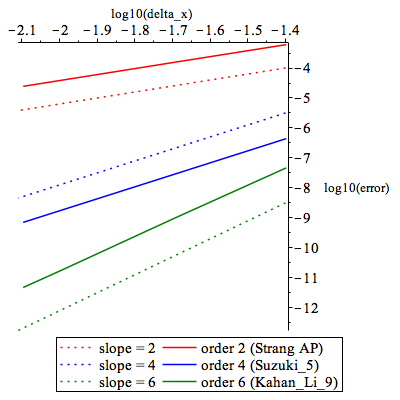}\includegraphics[width=6cm]{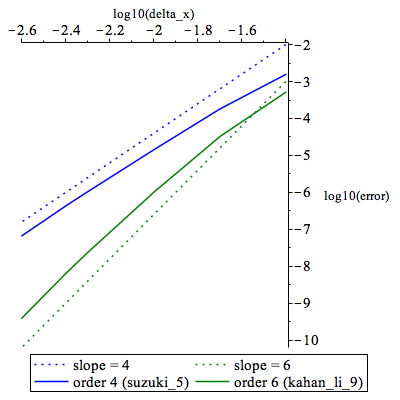}\caption{\label{fig:Convergence-study}Convergence study for several palindromic
methods, order 2 (red), 4 (blue) and 6 (green). The dotted lines are
reference lines with slopes 2, 4 and 6 respectively. Left: CFL number
$\beta=5.$ Right: CFL number $\beta=50$.}
\end{figure}

On Figure \ref{fig:Convergence-study} (left picture) we give the
results of the convergence study for the smooth solution. The considered
error is the $L^{2}$ norm of $\vf_{h}(\cdot,t_{\max})-\vf(\cdot,t_{\max}).$ 

We make the same experiment with $\beta=50$. The convergence study
for the Suzuki and Kahan-Li schemes is also presented on Figure \ref{fig:Convergence-study}
(right picture). At high CFL, not only the scheme remains stable,
but the high accuracy is also preserved.

\subsubsection{Behaviour for discontinuous solutions}

\begin{figure}
\centering{}\includegraphics[width=12cm]{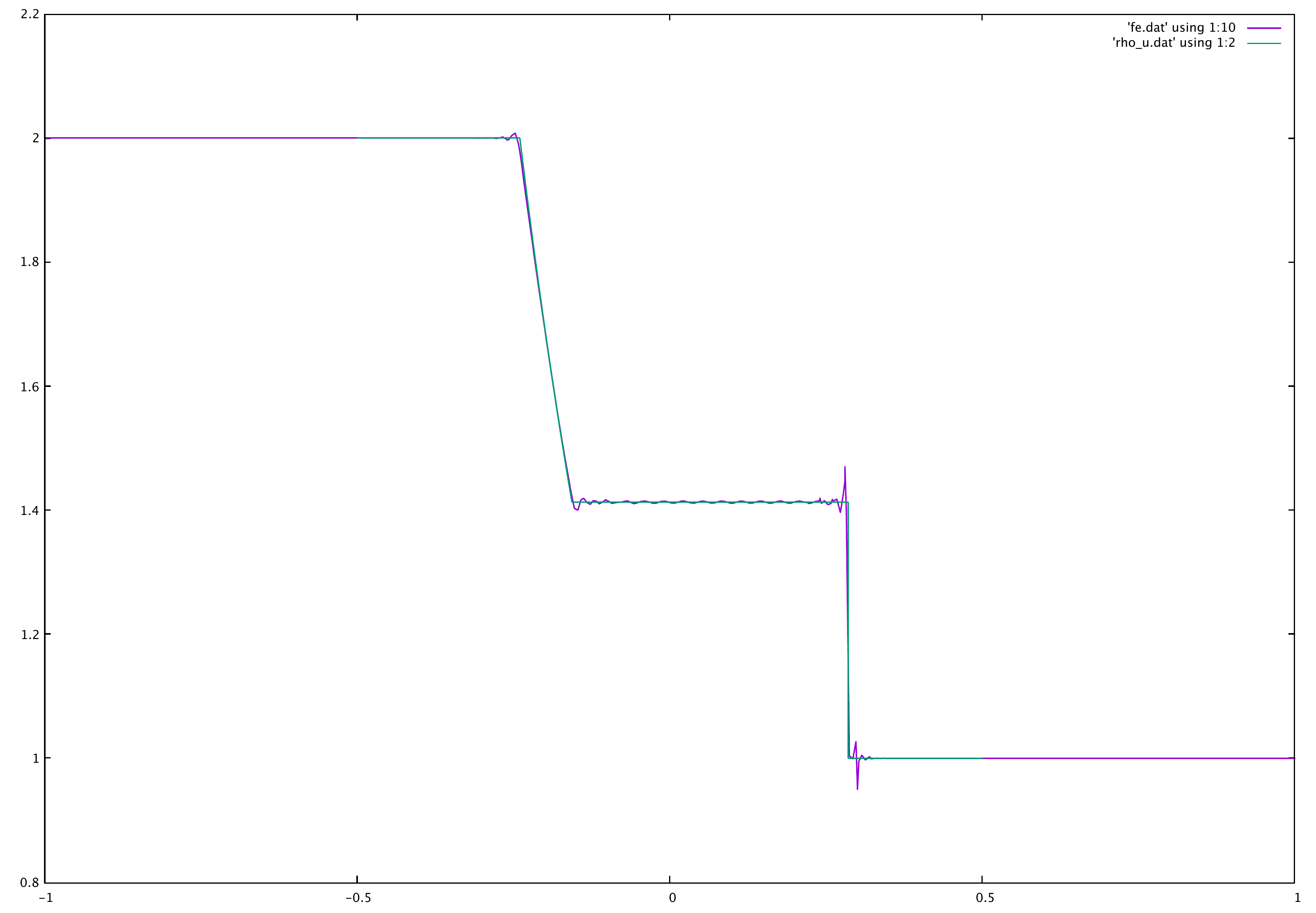}
\includegraphics[width=12cm]{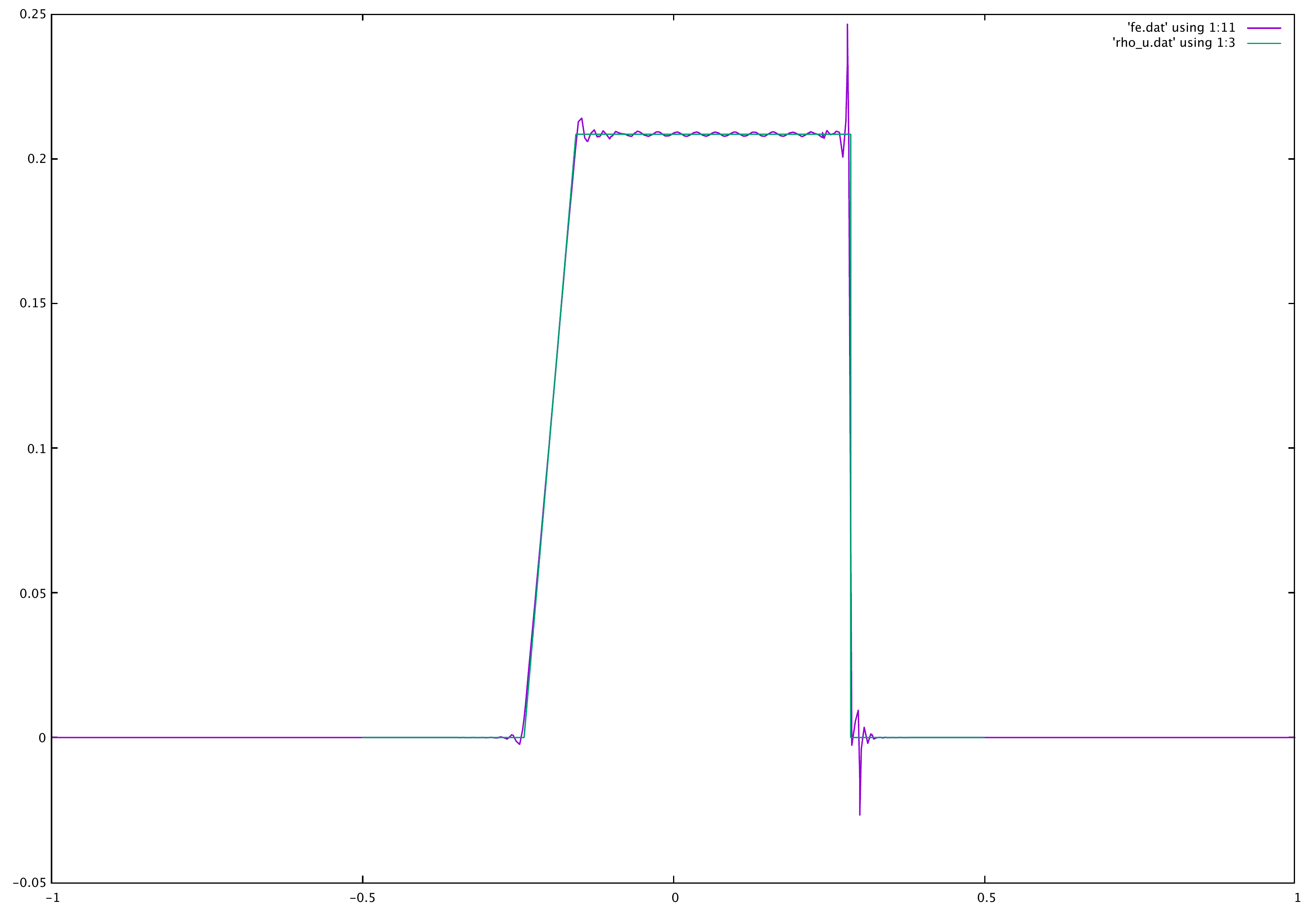}\caption{\label{fig:riemann_u}Riemann problem with $\tau=0$. Comparison of
the exact solution (green curve), and the numerical sixth-order solution
(purple curve). Left: density. Right:velocity.}
\end{figure}

We have also experimented the scheme for discontinuous solutions.
Of course, in this case the effective order of the method cannot be
higher than one and we expect Gibbs oscillations near the discontinuities.
On the interval $[a,b]=[-1,1],$ we consider a Riemann problem with
the following initial condition
\[
\rho(x,0)=\begin{cases}
2 & \text{if }x<0,\\
1 & \text{otherwise.}
\end{cases},\quad u(x,0)=0.
\]
We consider numerical results in the fluid limit $\tau=0$. On Figure
\ref{fig:riemann_u} we compare the sixth-order numerical solution
with the exact one at $t=t_{\max}=0.4$ for a CFL number $\beta=3$
and $N_{x}=100$ cells. We observe oscillations in the shock wave
and at the boundaries of the rarefaction wave, as expected. However,
we also observe that the high order scheme is able to capture a precise
rarefaction wave and the correct position of the shock wave. This
is a little bit surprising, because in presence of shock waves, the
Euler model is no more reversible and we solve it with a palindromic
time integrator method that has a reversible structure. The only dissipation
is provided by the upwind DG solver (see \eqref{eq:dissipation}). Apparently,
this slight dissipation is sufficient here for stabilizing the numerical
method.

\subsection{MHD flow}

We now consider a two-dimensional MagnetoHydroDynamics MHD model.
The $m=6$ unknowns of the model are the density $\rho(\vx,t)\in\mathbb{R}$,
the two-dimensional velocity vector $\v u(\vx,t)\in\mathbb{R}^{2}$, 
the two-dimensional magnetic field $\v B(\vx,t)\in\mathbb{R}^{2}$
and the total energy $Q(\vx,t)\in\mathbb{R}$. The pressure $p(\vx,t)\in\mathbb{R}$ is
given by a perfect gas pressure law
\[
p=(\gamma-1)\left(Q-\rho\frac{\v u\cdot\v u}{2}-\frac{\v B\cdot\v B}{2}\right),\quad\gamma=5/3.
\]

The conservative variables are
\[
\v w=(\rho,\rho\v u^{T},Q,\v B^{T})^{T}.
\]
For a two-dimensional direction vector $\v n=(n_{1},n_{2})^{T}$,
the MHD flux is then given by
\[
\v q(\v w)\v n=\left({\begin{array}{c}
{\rho\v u\cdot\v n}\\
{\rho(\v u\cdot\v n)\v u+(p+\frac{{\v B\cdot\v B}}{2})\v n-(\v B\cdot\v n)\v B}\\
{(Q+p+\frac{{\v B\cdot\v B}}{2})\v u\cdot\v n-(\v B\cdot\v u)(\v B\cdot\v n)}\\
{(\v u\cdot\v n)\v B-(\v B\cdot\v n)\v u}
\end{array}}\right).
\]
We approximate this hyperbolic system with the vectorial kinetic representation
described in Section \ref{subsec:vectorial-kinetic-scheme}.

The test-case is built upon a single vortex, which is a stationary
solution of the MHD system, to which a constant drift velocity is
added. In the moving frame centered on $\mathbf{r}_{O}(t)=t\mathbf{u}_{\text{\text{drift}}}$, with $\mathbf{u}_{\text{\text{drift}}} \in \mathbb{R}^{2}$,  the analytical solution reads in polar coordinates
\begin{align*}
\rho(r,\theta) & =\rho_{0},\\
\mathbf{u}(r,\theta) & =u_{0}[\mathbf{u}_{\text{\text{drift}}}+h(r)\mathbf{e}_{\theta}],\\
\mathbf{B}(r,\theta) & =b_{0}h(r)\mathbf{e}_{\theta},\\
p(r,\theta) & =p_{0}+\frac{b_{0}^{2}}{2}(1-h(r)),
\end{align*}
with $b_{0}=\rho_{0}u_{0}^{2}$. The results shown below are obtained with the parameter set $\rho_{0}=p_{0}=1,\ u_{0}=b_{0}=0.2,\ \mathbf{u}_{\text{drift}}=(1,1)^{T},\ h(r)=\exp[(1-r^{2})/2]$. The macromesh is the disk made from $20$
macrocells. Each macrocell is refined into $8\times8=64$ subcells
with fifth order basis functions, leading to $2304$ quadrature points
per macrocell. The minimal distance between two quadrature points
is around $h_{min}\approx0.02$. To each of the six scalar fields
$(\rho,\rho u_{x},\rho u_{y,}Q,B_{x}B_{y})$, we associate a four-velocity
$D2Q4$ model, with velocities $(-\lambda,0),(\lambda,0),(0,-\lambda),(0,\lambda)$,
so that there are $24$ kinetic fields. The velocity scale $\lambda$
is set to $4$. The vortex is initially centered in $\mathbf{r}_{O}=(0,0)^T$
at $t=0$ and we perform the simulation up to $t=1$. We test convergence
of the first, second and fourth order splitting schemes with time-steps
ranging from $\Delta t=0.2$ to $\Delta t=0.0125$. This leads to
kinetic CFL numbers ($\lambda\Delta t/h_{min}$) for the transport
subsets ranging from $40$ to $5$ for the first order splitting scheme,
$20$ to $2.5$ for the second order splitting scheme, and $13$ to
$1.6$ for the fourth order Suzuki scheme. 
\begin{figure}
\begin{centering}
\includegraphics[width=0.9\textwidth]{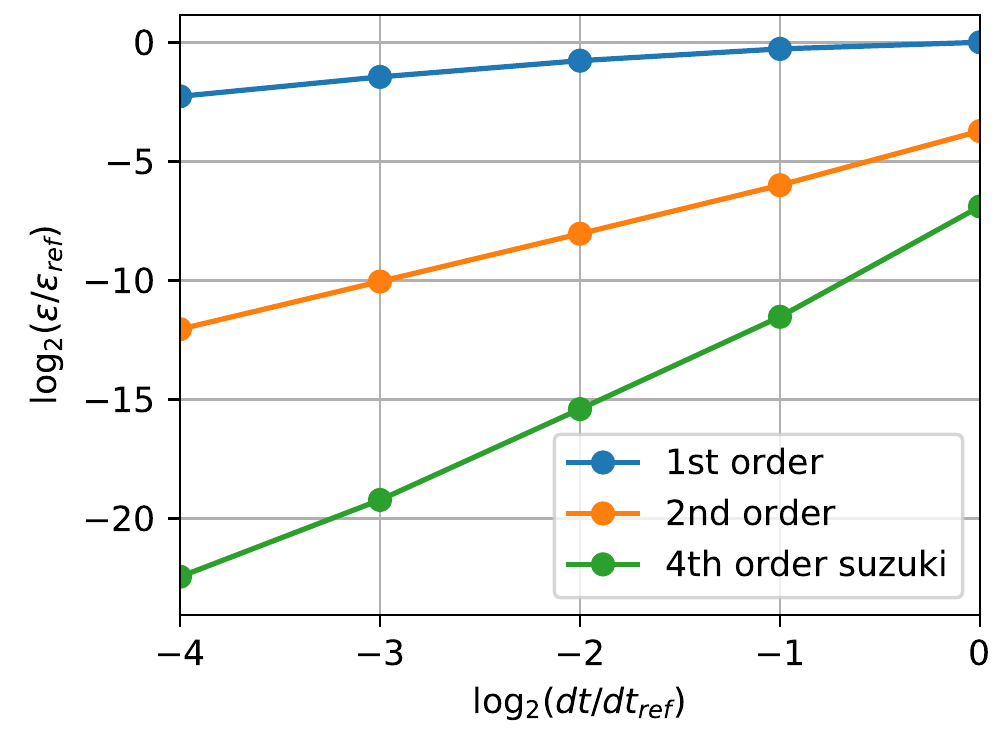}
\par\end{centering}
\caption{\label{fig:mjd_vortex_order}Drifting MHD vortex test-case. Convergence
of $L^{2}$error with respect to the analytical solution at $t=1$.
The reference error $\epsilon_{ref}$ of the log scale is the error
of the first order scheme for $\Delta t=0.2$.}
\end{figure}

\subsection{Flow past a cylinder (``thick'' boundary condition)\label{subsec:Flow-cylinder}}

We here consider the two-dimensional isothermal Euler equation and
its D2Q9 approximation presented in Example \ref{exa:2d_D2Q9} (see
Section \ref{subsec:LBM}).

\begin{figure}
\begin{centering}
\includegraphics[width=0.9\textwidth]{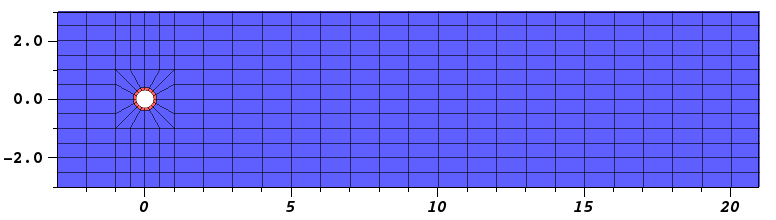}
\par\end{centering}
\caption{\label{fig:euler_cylinder_flow_mesh}Flow around a cylindrical obstacle.
Coarse macromesh with $328$ macrocells. The no-slip condition is
applied using a stiff relaxation in the thin annulus marked in red.}
\end{figure}

In this test case, we consider the flow of a fluid in a rectangular
duct with a cylindrical solid obstacle, as presented in Figure \ref{fig:euler_cylinder_flow_mesh}.
The simulation domain is the rectangle $[-3,21]\times[-6,6]$. The
cylindrical obstacle has radius $r=0.4$. At the boundary of the obstacle,
no-slip boundary conditions are applied. 

The initial condition is given by the constant state
\[
\rho(\vx,0)=1,\quad u(\vx,0)=u_{0}=0.05,\quad v(\vx,0)=0.
\]
 The boundary condition applied at the duct left entry ($x=-3$ axis)
for the whole simulation is 
\[
\rho(\vx,t)=1,\quad u(\vx,t)=\tanh(t/5)u_0,\quad v(\vx,t)=0.
\]
The progressive growth of the flow at the inlet mitigates the initial
unphysical transitory regime during which the initially uniform flow
adapts to the no-slip condition at the obstacle boundary. A second
consequence is the transition during the simulation between an initial
symmetric recirculation regime (with two vortices in the wake of the
obstacle) to the formation of von Karman streets \cite{GRUCELSKI2013406}.
\subsubsection{No-slip boundary condition\label{subsec:Thick-boundary-approach}}

To take into account the no-slip boundary condition around the obstacle,
we use a fictitious domain approach \cite{mittal2005immersed,angot2014optimal}. 

Our way to apply the fictitious domain approach is to first mesh a
thin annular shell (of width $0.1$), at the boundary of the obstacle.
The computational domain is thus enlarged with a small part of the
obstacle. The boundary condition is then applied by considering a
stiff penalization source term in the fluid equations
\[
\v s=-\kappa\left(\begin{array}{c}
0\\
\rho\v u
\end{array}\right),
\]
with $\kappa=0$ inside the fluid and $\kappa\gg1$ in the solid.
This amounts to considering the solid as a porous media with a very
small porosity. 

On the kinetic side, this source term can be represented in many different
ways. Our choice is to take
\[
\v g=-\kappa\left(0,f_{1}-f_{3},f_{2}-f_{4},f_{3}-f_{1},f_{4}-f_{2},f_{5}-f_{7},f_{6}-f_{8},f_{7}-f_{5},f_{8}-f_{6}\right).
\]
In other words, each component of the kinetic distribution associated
to a given lattice velocity relaxes toward the component associated
with the opposite velocity (see Figure \ref{fig:D2Q9_D3Q27_nodes}).

In practice, we observe a very fast decay of the velocity in the obstacle
as expected. In addition, this procedure is much more stable than
a Dirichlet type boundary condition (\ref{eq:boundary_conditions-1}).
See the numerical results presented below.

\subsubsection{Numerical simulations}

\begin{figure}
\begin{centering}
\includegraphics[width=0.9\textwidth]{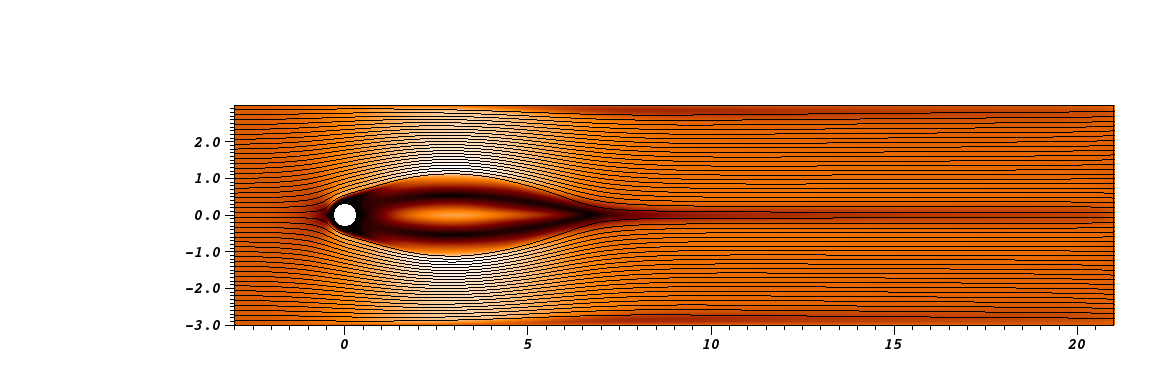}
\par\end{centering}
\begin{centering}
\includegraphics[width=0.9\textwidth]{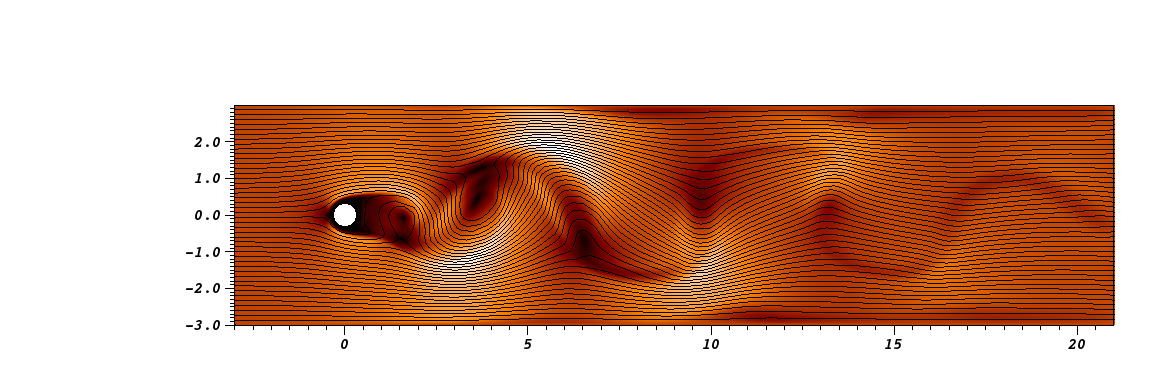}
\par\end{centering}
\begin{centering}
\includegraphics[width=0.9\textwidth]{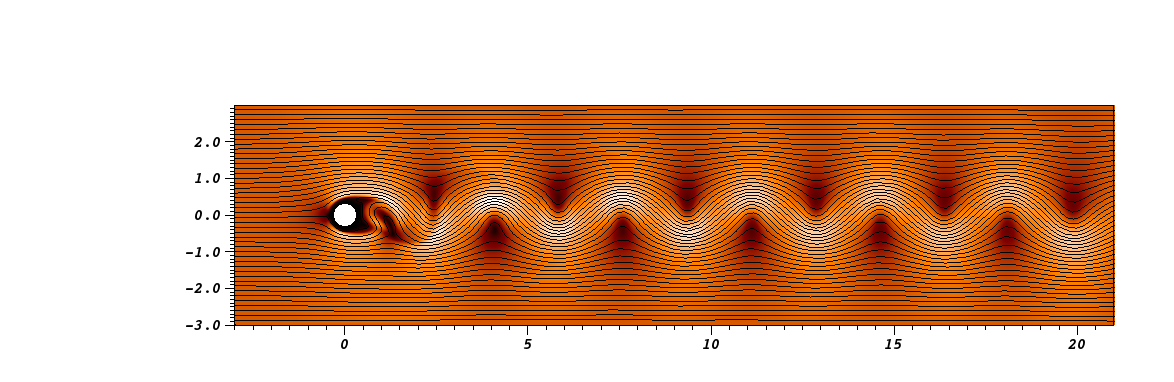}
\par\end{centering}
\caption{\label{fig:euler_cylinder_flow_stream_u}Flow around a cylindrical
obstacle. Velocity norm $\vert u\vert$ and streamlines at $t=50,t=120,t=340$. }
\end{figure}

\begin{figure}
\begin{centering}
\includegraphics[width=0.9\textwidth]{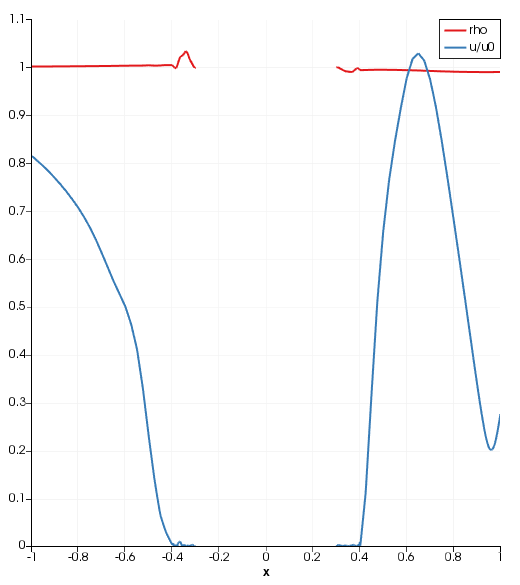}
\par\end{centering}
\caption{\label{fig:euler_cylinder_flow_urho_plot}Flow around a cylindrical
obstacle. Plot of $\vert u\vert/u_{0}$ and $\rho$ at $t=340$ on
the $x$ axis around the obstacle. }
\end{figure}

 The relaxation time has a finite but small value $\tau=0.0002$.
Accounting for the fact that for this model the dimensionless sound
speed is $c=1/\sqrt{3}$, the Mach number of the unperturbed flow is
approximately $ \frac{u_0}{c}\approx 0.087$. 
The simulation was performed on a macromesh
with $328$ macrocells; each macrocell contains $36\times36$ integration
points. The minimal distance between two integration points is about
$h_{min}\approx0.005$. The simulation was run with a time step $\Delta t=0.1$,
up to $t=3680$, about $7.5$ times the macroscopic transit time $L/u_{0}=480$.
For the $D2Q9$ kinetic model used herein, the maximal velocity modulus
is $\lambda_{max}=\sqrt{2}$. In the second order splitting scheme
the transport substep has $\text{\ensuremath{\Delta}t}=0.05$ at most,
so that the maximal kinetic transport CFL number is about $14$. On
Figure \ref{fig:euler_cylinder_flow_stream_u}, we show the streamlines
and velocity field norm at key points of the dynamics: at $t=50$
when the flow is still essentially symmetric, at $t=120$ after the
onset of the von Karman oscillations in the wake of the obstacle,
and at $t=340$ in the periodic oscillatory regime which starts at
about $t=220$. We observe that the geometry of the flow around the
obstacle is well preserved. The velocity is virtually null in the
thick boundary (Fig. \ref{fig:euler_cylinder_flow_urho_plot}). While
the density exhibits a small oscillation inside the thick boundary,
its value at the boundary of the computational domain is very close
to the nominal value. 

\subsection{Two-dimensional and three-dimensional two-fluid flow}

\subsubsection{Two-fluid flow with gravity}

We finally apply the methodology to a model of liquid-gas flow with
gravity. The model has been studied by several authors, see \cite{allaire2002five,chanteperdrix2002compressible,golay2007numerical}.
In dimension $2$ (resp. $3$), the $m=4$ (resp. $m=5$) unknowns
of the model are the mixture density $\rho(\vx,t)\in\mathbb{R}$,
the two-dimensional (resp. three-dimensional) velocity vector $\v u(\vx,t)\in\mathbb{R}^{2}$
(resp. $\mathbb{R}^{3}$) and the mass fraction of gas $\varphi(\vx,t)\in\mathbb{R}$.
The pressure of the mixture $p(\vx,t) \in \mathbb{R}$ is computed by
\[
p=p(\rho,\varphi)=\alpha\, p_{1}+(1-\alpha)\,p_{2},
\]
where $\alpha(\vx,t) \in[0,1]$ is the volume fraction of gas and $p_{1}(\vx,t)$, $p_{2}(\vx,t) \in  \mathbb{R}$
are the pressures of the gas and the liquid, respectively. For the
pressure of each fluid, we take
\[
p_{1}=p_{0}+c^{2}\left(\frac{\rho\varphi}{\alpha}-\rho_{0,1}\right),\quad p_{2}=p_{0}+c^{2}\left(\frac{\rho(1-\varphi)}{(1-\alpha)}-\rho_{0,2}\right).
\]
The physical constants of the models are the reference pressure $p_{0}$,
two reference densities for each fluid $\rho_{0,1}$ and $\rho_{0,2}$
and the sound speed $c$. Here the sound speed has no physical meaning.
It is chosen large enough in order that the flow can be considered
as almost incompressible. 

The volume fraction $\alpha$ is chosen in such a way that
\[
p_{1}=p_{2}.
\]
The conservative variables are
\[
\v w=(\rho,\rho\v u^{T},\rho\varphi)^{T}.
\]
The flux is given by
\[
\v q(\v w)\v n=\left(\begin{array}{c}
\rho\v u\cdot\v n\\
\rho(\v u\cdot\v n)\v u+p\v n\\
\rho\varphi u\cdot\v n
\end{array}\right),
\]
 and the source term
\[
\v s=\left(\begin{array}{c}
0\\
\rho\,\mathbf{g}\\
0
\end{array}\right).
\]

where $\mathbf{g}\in\mathbb{R}^{2}$ (resp. $\mathbb{R}^{3}$) is
the gravity vector in dimension $2$ (resp. $3$). In the following,
we consider the vectorial kinetic approximation as described in Section
\ref{subsec:vectorial-kinetic-scheme}. In dimension 2, unlike the
$D2Q9$ scheme, each of the four macroscopic fields is the sum of four
kinetic fields associated with the velocities $(-\lambda,0),(\lambda,0),(0,-\lambda)$ and
$(0,\lambda)$, leading to a total of $16$ kinetic fields. Similarly,
the three-dimensional kinetic relaxation model includes $30$ kinetic
fields. 

The macroscopic source term is represented by the kinetic source $\vg=\nabla_{\mathbf{w}}\vf^{\text{eq}}\v s$
of equation \eqref{eq:source_jacfeq_s}. 

In both test cases, we consider the growth of the Rayleigh-Taylor
instability: the light and heavy phases are initially well separated,
the heavy phase lying above the light one in the gravity field. At
the interface between the two phases, the mass fraction value drops
from $1$ to $0$ over the thin interface width. In order to avoid
Gibbs oscillations due to the large gradients, the relaxation time
is set to small (around $10^{-5})$ but finite values. The resulting
numerical viscosity smooths out the flow. 

\subsubsection{Two-dimensional Rayleigh-Taylor instability in an annulus.}

For the $2D$ model, we consider an annular domain of interior radius
$r_{min}=0.2$ and exterior radius $r_{max}=1$. The gravity field
is a radial one pointing inwards, i.e $\mathbf{g}=-g_{0}\mathbf{e}_{r}$
with $g_{0}=0.05$. The model parameters are 
\[
c=1,\quad p_{0}=1,\quad\rho_{0,1}=0.9,\quad\rho_{0,2}=1.1.
\]
In the initial unperturbed state, the heavy fluid lies ``above''
the light one in the annulus $r_{0}\le r\le1$, with $r_{0}=0.6.$
The interface is perturbed with a single azimuthal Fourier mode i.e.
\begin{equation}
r_{pert}=r_{0}+a\sin(m\theta).
\end{equation}
We will show here the results for a single $m=5$ azimuthal mode with
$a=0.01$. The sharpness of the transition from the light to the heavy fluid
is set using a hyperbolic tangent radial profile, so that the mass
fraction at $t=0$ reads
\begin{equation}
\varphi(r,\theta,t=0)=0.5\left(1-\tanh[(r-r_{pert}(\theta))/w_{pert}]\right).
\end{equation}
with $w_{pert}=0.02$.

The initial density is set so that each
of the pure phases is at mechanical equilibrium with the gravity field.
The macromesh is an annulus, discretized on a regular polar grid with
$n_{r}=5$ and $n_{\vartheta}=32$. Each macrocell is refined in $25\times10$
subcells with second order basis functions. With those parameters,
the minimal distance between two interpolation points is $h_{min}=0.002$.
The velocity scale parameter is set to $\lambda=2.5$. The time-step
is set to $\Delta t=0.01$. In the second order palindromic splitting
time scheme used here, the maximal time substep for the transport substep
is $0.5\Delta t=0.005$ so that the maximal $CFL$ number for
the transport of the kinetic fields is $6.25$. 

\begin{figure}
\begin{centering}
\includegraphics[width=0.5\textwidth]{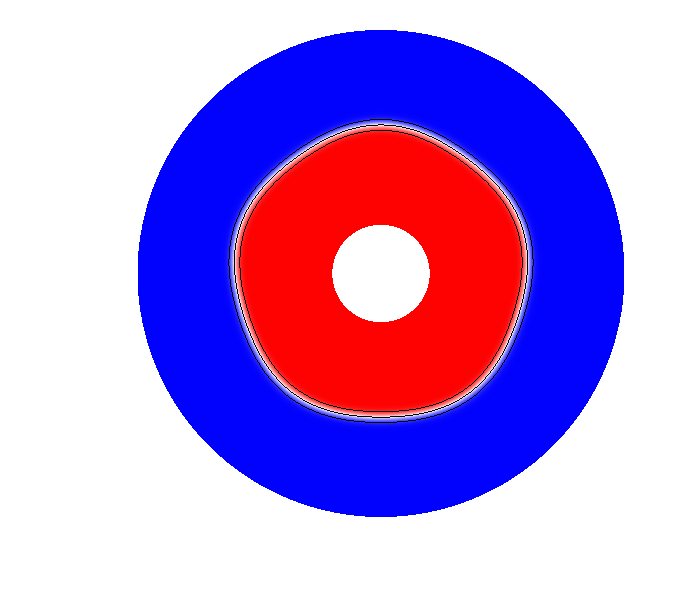}\includegraphics[width=0.5\textwidth]{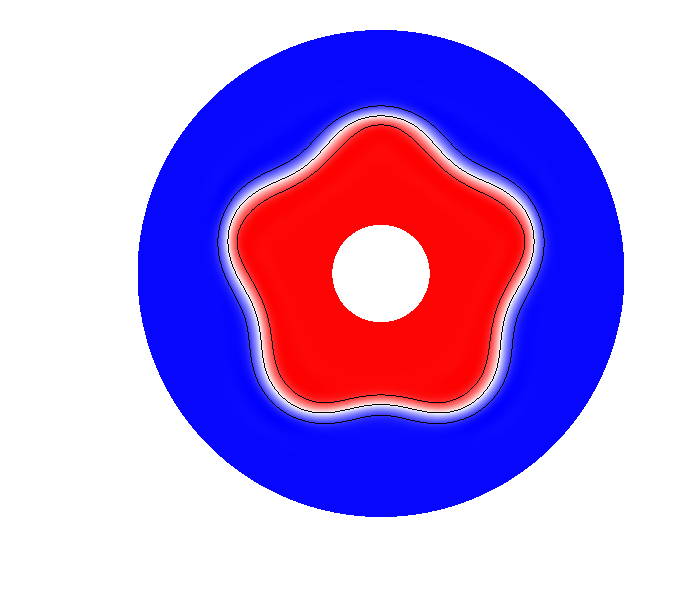}
\par\end{centering}
\begin{centering}
\includegraphics[width=0.5\textwidth]{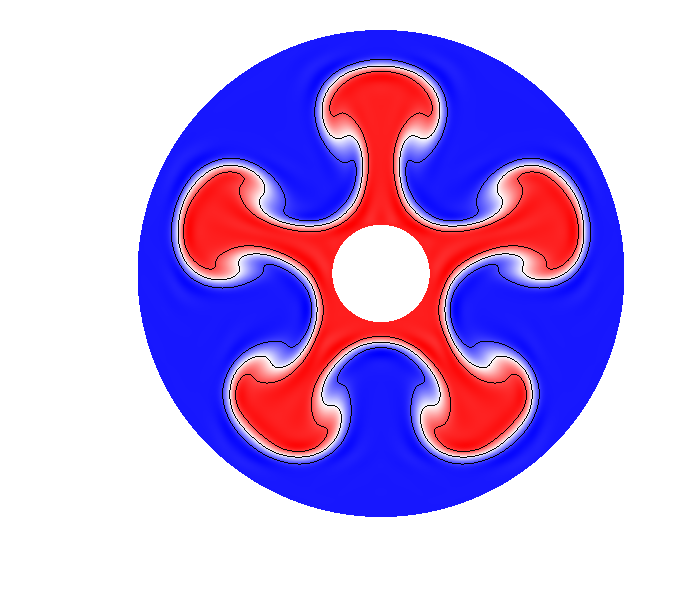}
\par\end{centering}
\caption{\label{fig:euler_mix_rt_disc}. Two-dimensional Rayleigh-Taylor instability
with a single $m=5$ azimuthal mode perturbation. Mass fraction $\varphi(\vx,t)$
at $t=[0,20,40]$. Black lines are iso-value contours at values $\varphi=0.1,0.5,0.9$.}
\end{figure}

\begin{figure}
\begin{centering}
\includegraphics[width=0.8\textwidth]{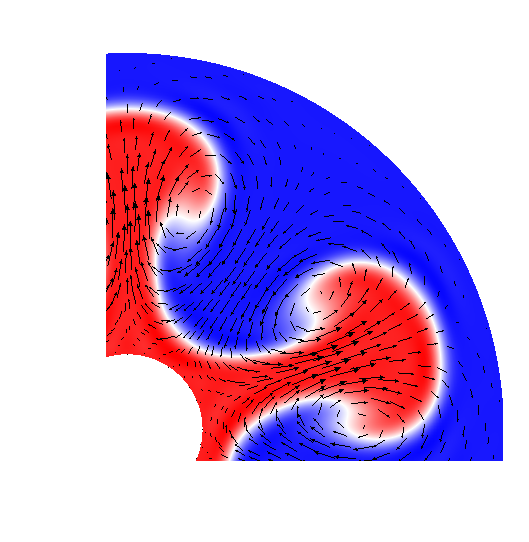}
\par\end{centering}
\caption{\label{fig:euler_mix_rt_disc2}Two-dimensional Rayleigh-Taylor instability
with a single $m=5$ azimuthal mode perturbation. Mass fraction $\varphi(\vx,t)$
and velocity field at $t=40$.}
\end{figure}

The evolution in time of the mass fraction $\varphi$ (see (Fig. \ref{fig:euler_mix_rt_disc}
and \ref{fig:euler_mix_rt_disc2}), the growth of the Rayleigh-Taylor
mushrooms is clearly visible. The dispersive errors entail a slight
excursion (a few percent) from the pure phase nominal values outside
of the interface zone. Those oscillations are mitigated by the diffusive
dissipation induced by the finite value of the relaxation time $\tau=0.00001$.
We do not concern ourselves here with the fine tweaking of parameters
or model required to strike a particular balance between the conservation
of the interface sharpness and the control of dispersive errors.

\begin{figure}
\begin{centering}
\includegraphics[width=0.5\textwidth]{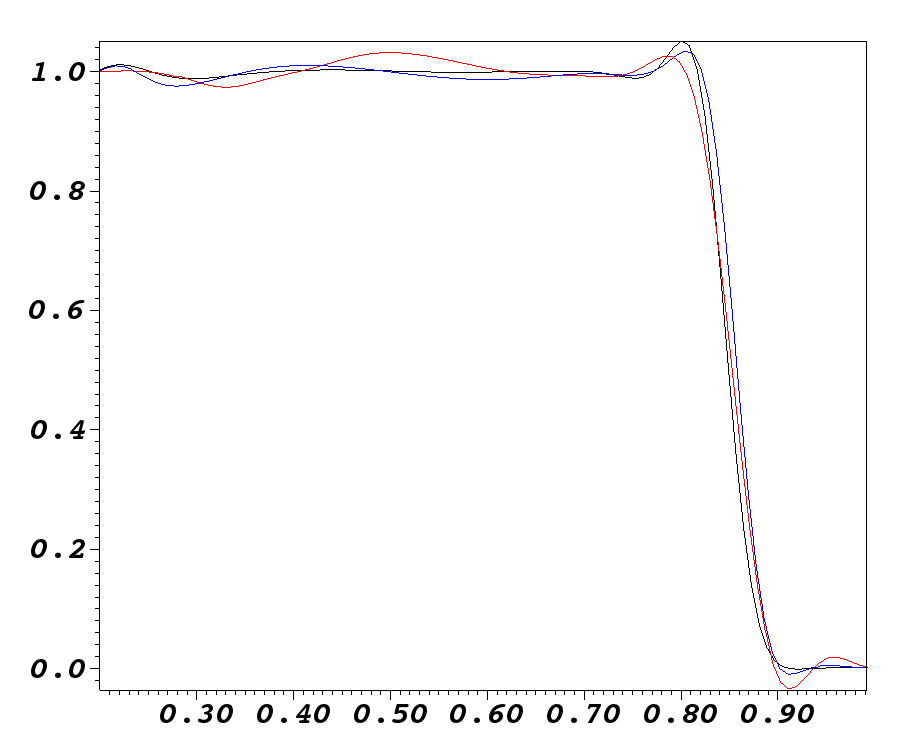}
\par\end{centering}
\caption{\label{fig:euler_mix_rt_disc3}Two-dimensional Rayleigh-Taylor instability
with a single $m=5$ azimuthal mode perturbation. Radial profiles of  mass fraction $\varphi(r,\theta,t)$
at time $t=40$ and azimuthal angles $\theta = \pi/2$ (continuous
black line), $\theta = \pi/2\pm2\pi/5$ (continuous and dotted blue lines), $\theta = \pi/2\pm4\pi/5$ (continuous
and dotted red lines). We observe small differences between the radial profiles,
because the kinetic velocity set is not aligned with the mesh, nor
with the mode.}
\end{figure}

Assuming a given overall accuracy, the question arises whether the
scarcity of the velocity set of the kinetic model induces geometrical
artifacts. For this particular test case the mesh, the $D2Q4$ velocity
set and the continuous system (perturbation included) are all symmetric
with respect to the $y$ axis. The fifth-order rotational symmetry
of the $m=5$ mode, is not preserved by the discrete velocity set
or the mesh though, and we can expect numerical anisotropy effects
to appear. On Figure \ref{fig:euler_mix_rt_disc3}, we compare radial
profiles of the mass fractions along the axes of the five Rayleigh-Taylor
mushrooms. The symmetry with respect to the vertical axis common to
both the excited mode and the velocity set is well preserved by the
scheme (the corresponding plots are indistinguishable), while the
$2\pi/5$ rotational symmetry is slightly broken due to the anisotropy
of the error. We conclude that the anisotropy of the kinetic velocity
set has only a very small effect.


\subsubsection{Three-dimensional Rayleigh-Taylor instability in a cylinder.}

\begin{figure}
\begin{centering}
\includegraphics[width=0.5\textwidth]{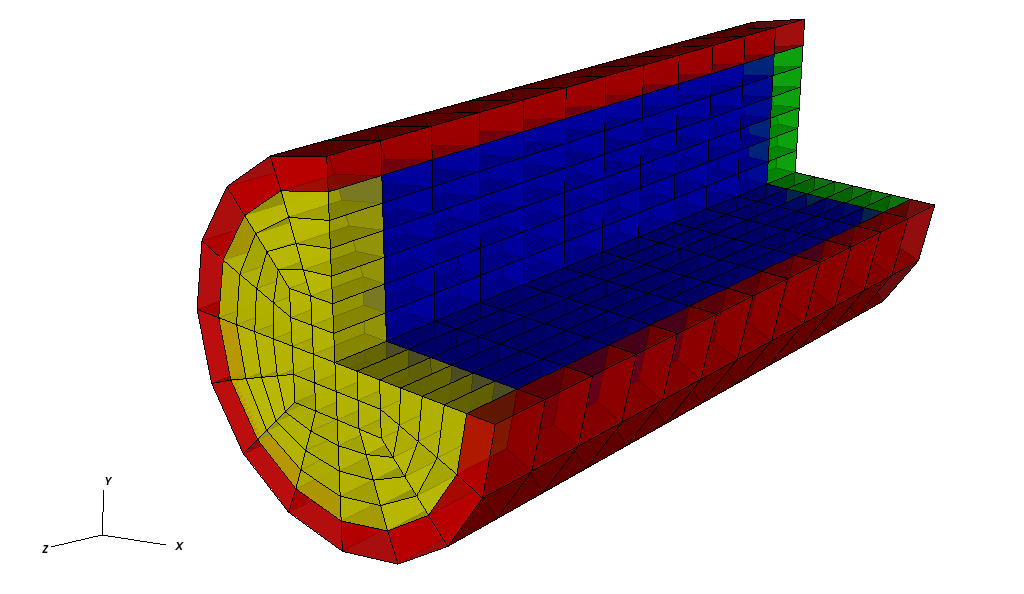}
\par\end{centering}
\caption{\label{fig:euler_mix_rt_cyl_mesh}Three-dimensional Rayleigh-Taylor
instability. Cylindrical macromesh with $1152$ macrocells. The physical
domain is the internal cylinder (blue internal zone and red and green
buffers). The annular buffer (red) is used to match the non-physical
imposed boundary condition on its exterior boundary at $r=0.2$ with
the physical boundary at $r=1$ using a volumic stiff relaxation operator. }
\end{figure}

We consider here the three-dimensional version of the two-fluids mixture
in a cylindrical duct of height $H=4.8$ and radius $R=1$ whose axis
is aligned with the constant gravity field $\mathbf{g}=-g_{0}\mathbf{e}_{z}$
with $g_{0}=0.04$. The coarse macromesh (Fig. \ref{fig:euler_mix_rt_cyl_mesh})
is composed of $1152$ macrocells, each of which is refined in $4\times4\times8=128$
subcells and second order basis functions. The buffer zones (yellow,
green and red on Fig. \ref{fig:euler_mix_rt_cyl_mesh}) are used to
apply boundary conditions on the internal cylindrical volume, with
a volumic relaxation operator. For this test case, the boundary conditions
at $z=\pm2.4$ are simply obtained by imposing the stationary equilibrium
state for each of the fluids, and no operator is required in the relevant
buffers (yellow and green on Fig. \ref{fig:euler_mix_rt_cyl_mesh}).
In order to mimic the effect of a solid duct at $r=1$, a no-slip
condition is applied on the horizontal directions and a slip condition
in the vertical direction.

The initial interface between the two fluids is the plane $h_{ref}=0.5$.
It is perturbed with a single bump centered in $r=0$ so that the
altitude $h_{pert}$ of the interface reads

\begin{equation}
h_{pert}(r)=h_{ref}-a\cos(2\pi r)\exp(-(r/d)^{2})
\end{equation}

with $h_{ref}=0.5,\ a=0.2,\ d=0.3$. The transition between the
light and heavy fluid is smoothed out using a hyperbolic tangent profile
of typical width $w_{pert}=0.05$, so that the mass fraction is given
in cylindrical coordinates by

\begin{equation}
\varphi(r,z,t=0)=0.5\left(1-\tanh(\frac{z-h_{pert}(r)}{w_{pert}})\right).
\end{equation}

\begin{figure}
\begin{centering}
\includegraphics[width=0.25\textwidth]{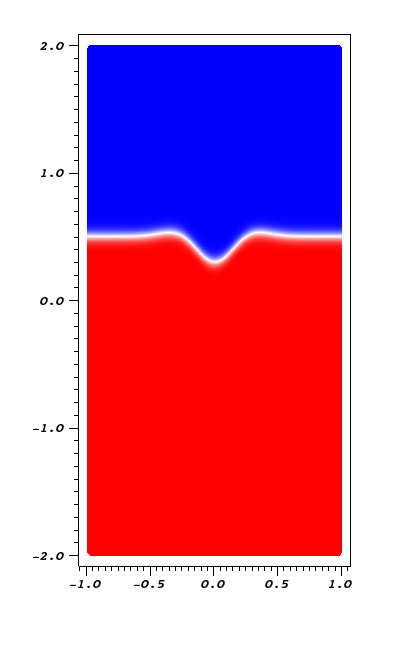}\includegraphics[width=0.25\textwidth]{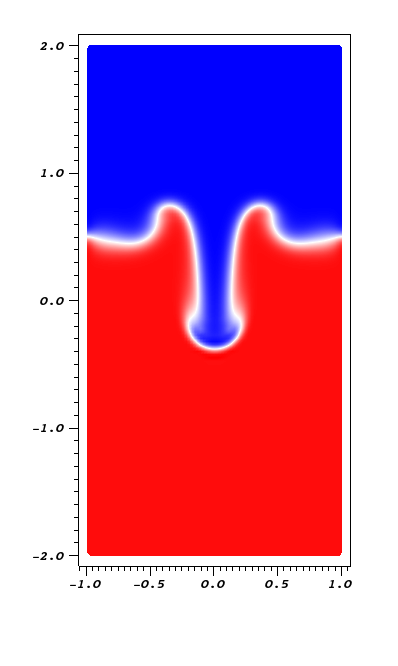}\includegraphics[width=0.25\textwidth]{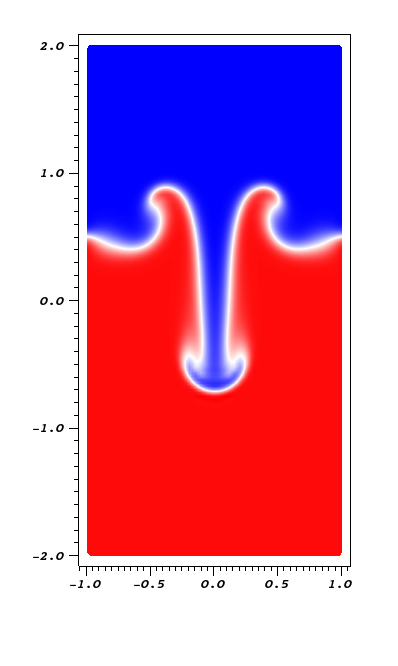}\includegraphics[width=0.25\textwidth]{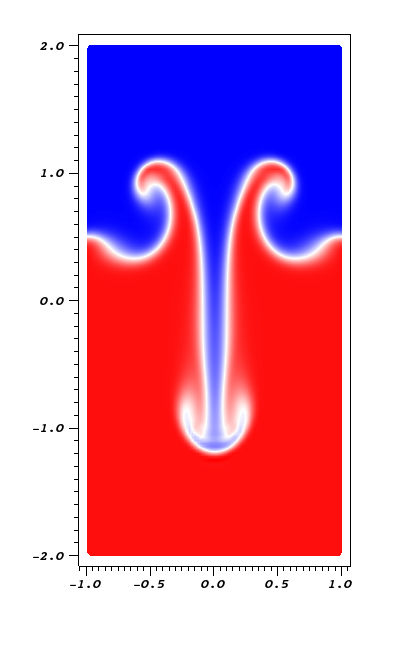}
\par\end{centering}
\caption{\label{fig:euler_mix_rt_cyl1}Three-dimensional Rayleigh-Taylor
instability. Planar cut at $y=0$ of the mass fraction $\phi$ at
$t=0$, $t=4$, $t=20$, $t=25$. }
\end{figure}

\begin{figure}
\begin{centering}
\includegraphics[width=0.25\textwidth]{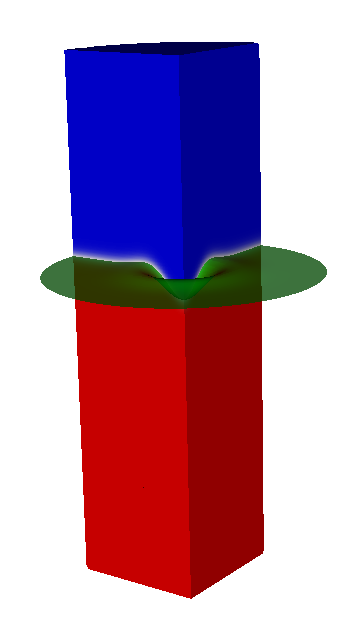}\includegraphics[width=0.25\textwidth]{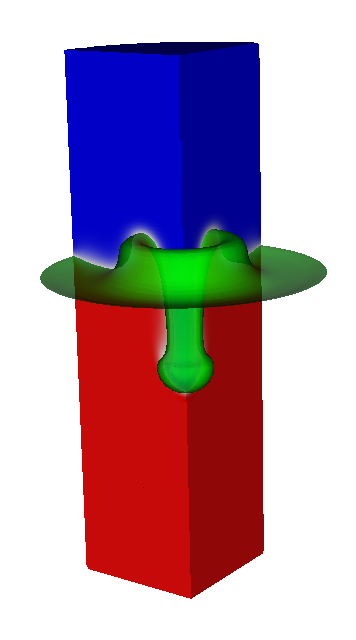}\includegraphics[width=0.25\textwidth]{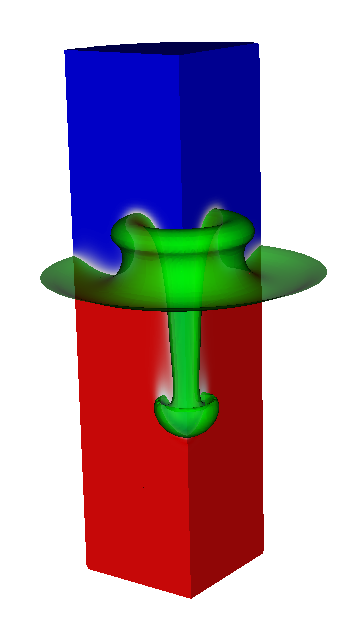}\includegraphics[width=0.25\textwidth]{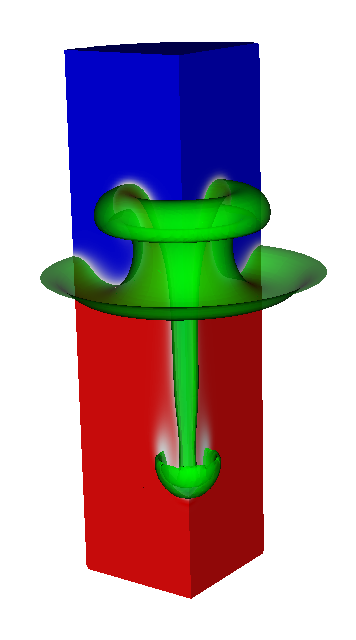}
\par\end{centering}
\caption{\label{fig:euler_mix_rt_cyl2}Three-dimensional Rayleigh-Taylor
instability. Pseudo-color slice and $\varphi=0.5$ iso-contour (green)
of the mass fraction $\varphi$ at times $t=0$, $t=4$, $t=20$, $t=25$. }
\end{figure}

On Fig. \ref{fig:euler_mix_rt_cyl1} and \ref{fig:euler_mix_rt_cyl2}, we
observe the development of the Rayleigh-Taylor mushroom. 

\section{Conclusion}

We have presented a new general implicit scheme, the Palindromic Discontinuous
Galerkin  (PDG) scheme, for solving systems of conservation laws. Despite
being formally implicit, the new scheme does not require costly linear
solver and has the complexity of an explicit scheme. We have also
proposed a new palindromic splitting algorithm that allows us to achieve
high accuracy in time, even in the stiff limit. The whole approach
remains stable and accurate, even at high CFL numbers.

We have validated the properties of the method on several one-dimensional
test cases. We have also tested the approach in higher dimensions,
and on different models of conservation laws coming from physics.

These first results are very promising for the future.

Many practical and theoretical works are still needed in order to
harness the full potential of the PDG method. The most important question
is to construct a methodology for handling general boundary conditions
in a stable way. A promising approach is to test and analyze in detail
the fictitious domain approach that we have sketched in this paper.
Another interesting direction of research would be to replace the
discontinuous Galerkin resolution of the transport equation by an
alternative method, such as semi-Lagrangian approaches. Finally, in
many applications, it is important to handle conservation laws with
small second order dissipative terms. Those dissipative terms can
be of a physical nature or serve a numerical purpose for avoiding
oscillations in shock waves, for instance. This can be achieved by
considering small, but non-vanishing, relaxation parameter $\tau>0$.
In this direction also, many useful practical extensions of the method
can be tested and analyzed.

\section{Appendix}

\subsection{Second order approximation\label{subsec:Second-order-approximation} }

For the sake of completeness, we recall the proof of (\ref{eq:macro_asymptic_first_order}).
\begin{proof}
We decompose $\vf$ into its equilibrium an non-equilibrium part setting
$\vf=\vf^{\text{eq}}+\vtf$. Substituting this formulation in the initial
kinetic system, and applying $P$ we get the equivalent coupled system 

\begin{equation}
\left\{ \begin{array}{cl}
{\displaystyle {\displaystyle }}\partial_{t}\vw+\Skd\partial_{k}P\vV^{k}\vf^{\text{eq}} & ={\displaystyle {\displaystyle }}-\partial_{k}P\vV^{k}\vtf+P\vg\\
{\displaystyle {\displaystyle }}\partial_{t}\vtf+\Skd\partial_{k}\vV^{k}\vtf & ={\displaystyle {\displaystyle }}-\tau^{-1}\vtf-\left[\partial_{t}\vf^{\text{eq}}+\Skd\partial_{k}\vV^{k}\vf^{\text{eq}}\right]+\vg
\end{array}\right..
\end{equation}
We now perform a formal expansion in $\tau$ of all quantities, with
$\vtf^{(0)}=0$: for instance, the kinetic source term is expanded
as $\vg=\vg^{(0)}+\tau\vg^{(1)}+\cdots$.

At the lowest order we have the limit system

\begin{equation}
\left\{ \begin{array}{ccc}
{\displaystyle {\displaystyle }}\partial_{t}\vw^{(0)}+\Skd\partial_{k}P\vV^{k}\vf^{\text{eq},(0)} & = & P\vg^{(0)}\\
\vtf^{(0)} & = & 0
\end{array}\right..
\end{equation}

At the first order we have

\[
\left\{ \begin{array}{ccc}
{\displaystyle {\displaystyle }}\partial_{t}\vw^{(1)}+\Skd\partial_{k}P\vV^{k}\vf^{\text{eq},(1)} & = & {\displaystyle {\displaystyle }}-\Skd\partial_{k}P\vV^{k}\vtf^{(1)}+P\vg^{(1)}\\
0 & = & {\displaystyle {\displaystyle }}-\vtf^{(1)}-\left[\partial_{t}\vf^{\text{eq},(0)}+\Skd\partial_{k}\vV^{k}\vf^{\text{eq},0}\right]+\vg^{(0)}
\end{array}\right.,
\]

and the second equation yields
\begin{equation}
\vtf^{(1)}={\displaystyle {\displaystyle }}-[\nabla_{\v w}\vf^{\text{eq},(0)}\partial_{t}\vw^{(0)}+\Skd\vV^{k}\nabla_{\v w}\vf^{\text{eq},(0)}\partial_{k}\vw^{(0)}]+\vg^{(0)}.
\end{equation}
Substituting the expression for $\partial_{t}\vw^{(0)}$ obtained
at the lowest order we obtain

\begin{multline}
\vtf^{(1)}=\Skd[\nabla_{\mathbf{w}}\vf^{\text{eq},(0)}P\vV^{k}\nabla_{\mathbf{w}}\vf^{\text{eq},(0)}-\vV^{k}\nabla_{\mathbf{w}}\vf^{\text{eq},(0)}]\partial_{k}\vw^{(0)}+\\{}
[\vg^{(0)}-\nabla_{\mathbf{w}}\vf^{\text{eq},(0)}P\vg^{(0)}].
\end{multline}
Recombining terms up to first order in $\tau$, and using the consistency
condition $P\vV^{k}\vf^{\text{eq}}=\vq^{k}$, the kinetic relaxation system
is consistent with

\begin{multline}
\partial_{t}\vw+\Skd\partial_{k}\vq^{k}(\vw)=\vs+\tau\Skd\Sjd\partial_{k}[\mathcal{D}^{kj}\partial_{j}\vw]+\\
\tau\Skd\partial_{k}P\vV^{k}[\nabla_{\mathbf{w}}\vf^{\text{eq}}\vs-\vg]+\mathcal{O}(\tau^{2}),
\end{multline}

with the diffusion tensor given by

\begin{equation}
\mathcal{D}^{kj}=P\vV^{k}\vV^{j}\nabla_{\mathbf{w}}\vf^{\text{eq}}-\nabla_{\mathbf{w}}\vq^{k}\nabla_{\mathbf{w}}\vq^{j}.
\end{equation}
\end{proof}

\section*{Bibliography}

\bibliographystyle{plain}
\bibliography{helluy-relax}

\end{document}